\documentclass[sort&compress,3p,fleqn]{elsarticle}
\usepackage{pdfpages,color,tikz,pdfsync,enumitem}
\usepackage{graphicx}
\usepackage{amsmath,amssymb,amsthm}
\usepackage[CJKbookmarks=True]{hyperref}

\allowdisplaybreaks[4]

\numberwithin{equation}{section}

\newtheorem{Def}{Definition}[section]
\newtheorem{Thm}[Def]{Theorem}
\newtheorem{Lem}[Def]{Lemma}
\newtheorem{Rmk}[Def]{Remark}

\newcommand{\R}{\mathbb{R}}
\newcommand{\J}{\mathbb{J}}
\newcommand{\mC}{\mathbb{C}}
\newcommand{\Z}{\mathbb{Z}}

\newcommand{\T}{\mathbb{T}^2}
\newcommand{\N}{\mathbb{N}}

\newcommand{\ud}{d}

\newcommand{\re}{\mathrm{Re}\,}
\newcommand{\im}{\mathrm{Im}\,}

\newcommand{\Hess}{\mathrm{Hess}}

\newcommand{\bvec}[1]{\boldsymbol{#1}}

\newcommand{\vj}{\bvec{j}}

\newcommand{\vq}{\bvec{q}}

\newcommand{\vx}{\bvec{x}}
\newcommand{\va}{\bvec{a}}

\newcommand{\val}{\bvec{a}^0}
\newcommand{\vb}{\bvec{b}}

\newcommand{\p}{\partial}

\newcommand{\Div}{\nabla\cdot}

\newcommand{\la}{\left\langle}
\newcommand{\ra}{\right\rangle}
\newcommand{\nn}{\nonumber\\}

\newcommand{\SE}{Schr\"{o}dinger equation}

\newcommand{\GL}{Ginzburg-Landau}
\newcommand{\wto}{\rightharpoonup}

\newcommand{\be}{\begin{equation}}
\newcommand{\ee}{\end{equation}}

\newcommand{\M}{W^{-1,1}(\T)}

\begin{document}

\begin{frontmatter}
\title{Quantized vortex dynamics of the complex \GL{} equation on torus}

\author[Zhu]{Yongxing Zhu\corref{cor}}
\ead{zhu-yx18@mails.tsinghua.edu.cn}

\address[Zhu]{Department of Mathematical Sciences, Tsinghua University, Beijing, 100084, China}

\cortext[cor]{Corresponding author.}

\begin{abstract}
We derive rigorously the reduced dynamical laws for quantized vortex dynamics of the  complex \GL{} equation on torus when the core size of vortex $\varepsilon\to 0$. The reduced dynamical laws of the complex \GL{} equation are governed by  a mixed flow of gradient flow and Hamiltonian flow which are both driven by a renormalized energy on torus. Finally, some first integrals and analytic solutions of the reduced dynamical laws are discussed.
\end{abstract}


\date{}
\end{frontmatter}
\section{Introduction}\label{sec:introduction}
We consider the vortex dynamics of the complex \GL{} equation on torus:
\begin{align}\label{eq:CGL2}
  \mu k_\varepsilon\p_tu^\varepsilon(\vx,t)+\lambda i\p_tu^\varepsilon(\vx,t)-\Delta u^\varepsilon(\vx,t)+\frac{1}{\varepsilon^2}(|u^\varepsilon(\vx,t)|^2-1)u^\varepsilon(\vx,t)=0,\vx\in\T,
\end{align}
with initial data
\begin{equation}\label{eq:initial CGL2}
  u^\varepsilon(\vx,0)=u^\varepsilon_0(\vx).
\end{equation}
Here $0<\varepsilon<<1$ is a parameter which is used to characterize the core size of quantized vortices,$k_\varepsilon=1/|\log \varepsilon|,\mu,\lambda$ are real parameters, $t$ is time, $\T=(\R/\Z)^2$ is the unit torus, $\vx=(x,y)^T\in\T$ is the spatial coordinate, $u^\varepsilon=u^\varepsilon(t)=u^\varepsilon(\vx,t)$ is a complex-valued wave function or order parameter and $u^\varepsilon_0:=u^\varepsilon_0(\vx,t)$ is complex-valued.


We define the \GL{} functional (energy) by \cite{BaoTang2014NumericalStudyQuantizedVorticesSchrodinger,JerrardSoner1998DynamicsGinzburgLandauVortices,Lin1996DynamicGinzburgLandau}
\begin{equation}\label{eq:def of E2}
  E^\varepsilon(u^\varepsilon(t)):=\int_{\T}e^\varepsilon(u^\varepsilon(\vx,t))d\vx,
\end{equation}
the momentum by \cite{CollianderJerrard1999GLvorticesSE}
\begin{equation}\label{eq:def of Q2}
  \bvec{Q}(u^\varepsilon(t)):=\int_{\T}\vj(u^\varepsilon(\vx,t))d\vx,
\end{equation}
where energy density $e^\varepsilon(v)$, current $\vj(v)$ and Jacobian $J(v)$ is defined as follows: for any complex-valued function $v:\T\to\mC^2$,
\begin{equation}\label{eq:def of e j J}
  \vj(v):=\im (\overline{v}\nabla v), \quad
e(v):=\frac{1}{2}|\nabla v|^2+\frac{1}{4\varepsilon^2}(1-|v|^2)^2,\quad J(v)=\frac{1}{2}\nabla\cdot(\J \vj(v))=\im (\partial_x\overline{v}\,\partial_y v),
\end{equation}
with $\overline{v},\im v$ denoting the complex conjugate and imaginary part of $v$, respectively, and 
\begin{equation}
  \J=\left(\begin{array}{cc}0&1\\-1&0 \end{array}\right).
\end{equation}

\eqref{eq:CGL2}  describes weakly nonlinear spatiotemporal phenomena in extended continuous media as an amplitude equation in physics \cite{AransonKramer2002ComplexGinzburgLandau,HohenbergKrekhov2015CGLinPhaseTransition}. And \eqref{eq:CGL2} is also used to model the competition of many species in some ecological systems \cite{MowlaeiRomanPleimling2014CGLinCompetitionSpecies}. When $\mu=0$ and $\lambda=1$, \eqref{eq:CGL2} collapses to the well-known nonlinear \SE{}, also known as Gross-Pitaevskii equation, which is used to describe superfluidity and Bose-Einstein condensate \cite{BaoBuZhang,Anderson2010ExperimentsVorticeSuperfluid,LinXin1999,Pitaevskii1961,Serfaty2017}:
\begin{align}\label{eq:NLS}
  i\p_tu^\varepsilon(\vx,t)-\Delta u^\varepsilon(\vx,t)+\frac{1}{\varepsilon^2}(|u^\varepsilon(\vx,t)|^2-1)u^\varepsilon(\vx,t)=0,\vx\in\T.
\end{align}
And when $\lambda=0,\mu=1$, \eqref{eq:CGL2} collapses to the \GL{} equation, which is used to describe type-II superconductivity \cite{Anderson2010ExperimentsVorticeSuperfluid,BaoTang2013NumericalQuantizedVortexGinzburgLandau,JerrardSoner1998DynamicsGinzburgLandauVortices,Lin1996DynamicGinzburgLandau}:
\begin{align}\label{eq:GL2}
  k_\varepsilon\p_tu^\varepsilon(\vx,t)-\Delta u^\varepsilon(\vx,t)+\frac{1}{\varepsilon^2}(|u^\varepsilon(\vx,t)|^2-1)u^\varepsilon(\vx,t)=0,\vx\in\T.
\end{align}

An important phenomenon of superfluid, type-II superconductivity, and other physical system related to \eqref{eq:CGL2} is the appearance of vortices, which are topological defects \cite{AransonKramer2002ComplexGinzburgLandau,Anderson2010ExperimentsVorticeSuperfluid}. In two-dimensional cases, quantized vortices are the zeros of order parameters with nonzero winding numbers (degrees). The appearance of quantized vortices is widely observed in  experiments related to superfluid, type-II superconductivity, and other physical system related to \eqref{eq:CGL2} \cite{BewleyLathropSreenivesan2006VisuaQuantizedvortex,Anderson2010ExperimentsVorticeSuperfluid}. From experiments and analytic studies \cite{Mironescu1995StabilityRadialSolutionGinzburgLandau,Anderson2010ExperimentsVorticeSuperfluid}, the vortex is stable only when its degree $d=\pm 1$.

The interaction and reduced dynamical laws of quantized vortices of \eqref{eq:CGL2} (and \eqref{eq:NLS}, \eqref{eq:GL2}) are widely studied in the past decades. For the results in the whole space $\R^2$ or on a bounded domain with Neumann boundary conditions or Dirichlet boundary conditions or the sphere $\mathbb{S}^2$, one can see \cite{test,Rubinstein1991SelfinducedmotionLinedefects,Neu1990VorticesComplexScalarfield,E1994DynamicsGinzburgLandau,BaoTang2013NumericalQuantizedVortexGinzburgLandau,BaoTang2014NumericalStudyQuantizedVorticesSchrodinger,SandierSerfaty2004EstimateGinzburgLandau,Lin1996DynamicGinzburgLandau,BetheulOrlandiSmets2006GinzburgLandautoMotionByMeanCurvature,ChenSternberg2014VortexDynamicsManifold,BaoBuZhang,Bethuel2008,DuJuGInzburgLandauSphereNumerical,Lin1999,LinXin1999,Zhang2007,Miot2009VorticesComplexGinzburgLandau,KurzkeMelcherMoserSpirn2009VortexDynamicsGinzburgLandauMixedFlow,AguarelesChapmanWitelski2008INteractionSprialCGL,AguarelesChapmanWitelski2020SprialWavesCGL,JianLiu2006GinzburgLandauVortexMeanCurvatureFlow,Jerrard1999GinzburgLandauwave,JerrardSoner1998DynamicsGinzburgLandauVortices,JerrardSpirn2008RefinedEstimateGrossPitaevskiiVortices,SandierSerfaty2004GammaConvergenceGinzburgLandau} and references therein for details. It was shown in the literature that for \eqref{eq:CGL2}  on $\Omega$, which was taken as $\R^2$ \cite{Miot2009VorticesComplexGinzburgLandau} or a bounded domain \cite{KurzkeMelcherMoserSpirn2009VortexDynamicsGinzburgLandauMixedFlow}, assuming that initial data $u_0^\varepsilon$ possesses $M$ vortices $\va_1^{0,\varepsilon},\cdots,\va_M^{0,\varepsilon}\in\Omega$ with degrees $d_1,\cdots,d_M\in\{\pm 1\}$ satisfying
satisfies
\begin{equation}
\va_j^{0,\varepsilon}\ne\va_k^{0,\varepsilon},\quad 1\le j\ne k\le M
\end{equation}
with 
\begin{equation}
  \va_j^0=\lim_{\varepsilon\to 0}\va_j^{0,\varepsilon},\quad \text{and}\ \va_j^0\ne \va_k^0,\ 1\le j\ne k\le M,
\end{equation}
then before the collision occurs (Collision occurs commonly as shown in \cite{BSX} and Subsection \ref{sec:solutions} if $\mu\ne 0$ and the degrees of vortices are not the same.), $u^\varepsilon(t)$ also possess $M$ vortices $\va_1(t),\cdots,\va_M(t)$ with degrees $d_1,\cdots,d_M$ satisfying
\begin{equation}
  \va^\varepsilon(t):=(\va_1^\varepsilon(t),\cdots,\va_M^\varepsilon(t))^T\to \va=\va(t)=(\va_1(t),\cdots,\va_M(t))^T\ \text{as} \ \varepsilon\to 0,
\end{equation}
where $\va$ satisfies 
\begin{equation}
  \va_j(t)\ne \va_k(t),\quad 1\le j\ne k\le M,
\end{equation}
and the reduced dynamical laws:
\begin{equation}
  \mu\dot{\va_j}-d_j\lambda\J\dot{\va_j}=-\frac{1}{\pi}\nabla_{\va_j}W_\Omega(\va).
\end{equation}
Here, $W_\Omega$ is the renormalized energy on $\Omega$ of the following form:
\begin{equation}\label{eq:RDL on Omega}
  W_\Omega(\va)=-\pi\sum_{1\le j\ne k\le M}d_jd_k\log|\va_j-\va_k|+\text{Reminding term determined by} \ \Omega.
\end{equation}

 \cite{CollianderJerrard1999GLvorticesSE,ZhuBaoJian} analytically studied the quantized vortex dynamics of \eqref{eq:NLS} on torus. These two papers showed two differences between the vortex dynamics on torus and on simply connected $\Omega$:
\begin{enumerate}[label=(\roman*)]
  \item Since $\T$ is a compact manifold, the summation of all degrees must be $0$ \cite{CollianderJerrard1999GLvorticesSE}. Noting that the degrees of vortices must be $\pm 1$ for stability, we can assume the number of vortices $M=2N(N\in\N)$ and
  \begin{equation}\label{eq:degrees2}
      d_1=\cdots= d_n=1,\quad d_{N+1}=\cdots=d_{2N}=-1.
  \end{equation}
  \item The reduced dynamical laws are determined not only by the position of vortices as in \eqref{eq:RDL on Omega}, but also by:
  \begin{equation}
    \bvec{Q}_0:=\lim_{\varepsilon\to 0}\bvec{Q}(u^\varepsilon_0).
  \end{equation}
\end{enumerate}

The main purpose of this paper is to give the dynamical law of vortices for \eqref{eq:CGL2} as $\varepsilon\to 0$ involving the influence of the $\bvec{Q}_0$.

To illustrate our result, we introduce some notations.

We first define
\begin{equation}
  (\T)^{2N}_*=\{(\vx_1,\cdots,\vx_{2N})\in(\T)^{2N}:\vx_j\ne\vx_k\ \text{for} \ 1\le j<k\le 2N\}.
\end{equation}
Then for any $\va=(\va_1,\cdots,\va_{2N})^T\in(\T)_*^{2N}$ and $\bvec{q}\in2\pi\sum_{j=1}^{2N}d_j\va_j+2\pi\Z^2$, the renormalized energy is defined by \cite{IgnatJerrard2021Renormalizedenergymanifold}
\begin{equation}\label{eq:def of W}
  W(\va;\bvec{q})=-\pi\sum_{1\le j\ne k\le 2N}d_jd_kF(\va_j-\va_k)+\frac{1}{2}|\bvec{q}|^2,
\end{equation}
with $F(\vx)$ the solution of
\begin{equation}\label{eq:define of F}
\Delta F(\vx)=2\pi (\delta(\vx)-1),\ \vx\in \T, \quad {\rm with}\quad \int_{\T}F(\vx)d\vx=0.
\end{equation}
Define 
\begin{equation}
  \gamma=\lim_{\varepsilon\to 0}\left(\inf_{v\in H^1_g(B_1(\bvec{0}))}\int_{B_1(\bvec{0})}\left[\frac{1}{2}|\nabla v|^2+\frac{1}{4 \varepsilon^2}(|v|^2-1)^2\right]d\vx-\pi\log\frac{1}{\varepsilon}\right),
\end{equation}
where the function space $H^1_g(B_1(\bvec{0}))$ is defined as
\[
  H_g^1(B_1(\bvec{0}))=\left\{v\in H^1(B_1(\bvec{0}))\left| v(\vx)=g(\vx)=\frac{x+iy}{|\vx|}\ \text{for}\  \vx\in \p B_1(\bvec{0})\right.\right\}.
\]
Then we introduce 
\begin{equation}\label{eq:def of Wepsilon}
  W_\varepsilon(\va;\bvec{q})=2N\left(\pi\log\frac{1}{\varepsilon}+\gamma \right)+W(\va;\bvec{q}).
\end{equation}

For simplicity, we will concentrate on the case
\begin{equation}\label{eq:mu=1}
  \mu=1
\end{equation}
in the rest of this paper.

Our main result is the following theorem:
\begin{Thm}[Reduced dynamical laws of the complex \GL{} equation]\label{thm:dynamic CGL2}
Assume there exists $\va^0=(\va_1^0,\dots,$ $\va^0_{2N})^T\in(\T)^{2N}_*$, $\bvec{q}_0\in 2\pi\sum_{j=1}^{2N}d_j\va_j^0+2\pi\Z^2$ such that the initial data $u^\varepsilon_0$ of \eqref{eq:CGL2} satisfies \eqref{eq:mu=1} and 
 \begin{equation}\label{con:limit of integer of J(varphi)}
    J(u^\varepsilon_0)\to\pi\sum_{j=1}^{2N}d_j\delta_{\va^0_j}\ \text{in}\ \M:=(C^1(\T))',\quad \lim_{\varepsilon\to 0}(E^\varepsilon(u^\varepsilon_0)- W_\varepsilon(\val;\vq_0))=0,\quad \bvec{Q}_0= \J \bvec{q}_0.
 \end{equation}
Then there exist Lipschitz paths $\va_j:[0,T)\to\T$, such that
\[
    J(u^\varepsilon(\vx,t))\to \pi\sum_{j=1}^{2N}d_j\delta_{\va_j(t)}\ \text{in}\ \M,\quad k_\varepsilon e^\varepsilon(u^\varepsilon(\vx,t))\to\pi\sum_{j=1}^{2N}\delta_{\va_j(t)}\ \text{in}\ \M,
\]
 and $\va=\va(t)=(\va_1(t),\cdots,\va_{2N}(t))^T$ satisfies 
\begin{align}
  \dot{\va}_j-\lambda d_j\J\dot{\va}_j=-\frac{1}{\pi}\nabla_{\va_j}W(\va;\vq(\va))\label{eq:CGLODE}
\end{align}
with initial data
\begin{equation}\label{eq:initial data CGL}
  \va_j(0)=\va_j^0,
\end{equation}
where $\vq(\va)=\vq(\va(t))$ is the lift of $\Pi(2\pi\sum_{j=1}^{2N}d_j\va_j(t))$, i.e. the unique continuous map $\bvec{q}(\va(t)):\R\to \R^2$ satisfying $\vq(\va(0))=\vq_0$ such that the following diagram commute.
\begin{center}
\begin{tikzpicture}
\draw[->](0,0)--(3.7,0);
\draw[->](0.2,-0.1)--(3.8,-1.9);
\draw[->](4,-0.3)--(4,-1.7);
\fill[white](0,0)circle(0.3);
\node at(0,0) {$\R$};
\node at(4.1,0) {$\R^2$};
\node at(4.1,-1.9) {$\T$};
\node at(2,0.2){$t\mapsto \bvec{q}(\va(t))$};
\node at(0,-1){$t\mapsto \Pi(2\pi\sum_{j=1}^{2N}d_j\va_j(t))$};
\node at(4.2,-1){$\Pi$};
\end{tikzpicture}
\end{center}
Here, $\Pi$ is  the canonical projection from $\R^2$ to $\T$.
\end{Thm}
We remark here that the constraint of $\vq_0$ is natural, since \eqref{con:limit of integer of J(varphi)} and Lemma \ref{lem:limitation of limit of current} imply that 
\begin{equation}\label{eq:limit of momentum of u_02}
\bvec{Q}_0=\J\vq_0\in2\pi\J\sum_{j=1}^{2N}d_j\va_j^0+2\pi\Z^2.
\end{equation}

The paper is organized as follows. In Section \ref{sec:Preliminaries2}, we will study the $\mathbb{S}^1$-valued function on torus and give some properties of the canonical harmonic map and the renormalized energy on torus. In Section \ref{sec:CGL}, we prove the existence of vortex paths of the solution of \eqref{eq:CGL2} and the convergence of the current of the solution. In Section \ref{sec:dynamics}, we establish the reduced dynamical laws \eqref{eq:CGLODE}. In Section \ref{sec:properties of CGLODE}, we give some first integrals of \eqref{eq:CGLODE} and some analytic solutions of \eqref{eq:CGLODE} with initial data with special symmetry.

\section{Preliminaries}\label{sec:Preliminaries2}

\subsection[S1-valued function on torus]{\texorpdfstring{$\mathbb{S}^1$}{}-valued function on torus}\label{sec:S1-valued}

\begin{Lem}\label{eq:momentum of S1function}
If $\mathbb{S}^1$-valued function $u\in H^1_{loc}(\T\setminus\{\va_1,\cdots,\va_{2N}\})\cap W^{1,1}(\T)$ satisfies
\begin{equation}\label{eq:def of Jacpobian}
  J(u)=\pi\sum_{j=1}^{2N}d_j\delta_{\va_j}\ \text{in}\ \M,
\end{equation}
where $\va=(\va_1,\cdots,\va_{2N})\in(\T)_*^{2N}$,  then we have
\begin{equation}\label{eq:limitation of momentum}
  \int_{\T}\vj(u)d\vx\in2\pi\J\sum_{j=1}^{2N}d_j\va_j+2\pi\Z^2.
\end{equation}

\end{Lem}

\begin{proof}
Without loss of generality, we can assume $\va_j=(x_j,y_j)\in(0,1)^2$. Since $u$ is an $\mathbb{S}^1$-valued function,  we can denote
\begin{equation}
  u(\vx)=e^{i\Theta(\vx)},
\end{equation}
which together with \eqref{eq:def of e j J} implies
\begin{equation}\label{eq:current of eta}
  \vj(u)=\im(\overline{u}\nabla u)=\nabla \Theta.
\end{equation}

Since $u$ is a function on torus, \eqref{eq:current of eta} implies that there exists a constant $k_1\in\Z$ such that
\begin{equation}\label{eq:int j on [0,1]xx}
  \int_0^1j_1(u(x,0))dx=\int_0^1\p_x\Theta(x,0)dx=\Theta(1,0)-\Theta(0,0)= 2\pi k_1.
\end{equation}
Then, via \eqref{eq:def of Jacpobian} and the definition of Jacobian \eqref{eq:def of e j J}, we have for any $y\in[0,1]$ and $\Gamma=\mathbb{S}^1\times[0,y]$,
\begin{align}\label{eq:difference of int j}
  \int_0^1j_1(u(x,0))dx-&\int_0^1j_1(u(x,y))dx=\int_{\p \Gamma}\vj(u(\vx))\cdot \bvec{\tau}ds\nn 
  =&\int_{\Gamma}\nabla\cdot(\J\vj(u(\vx)))d\vx=2\pi\int_{\Gamma}\sum_{j=1}^{2N}d_j\delta_{\va_j}d\vx\nn
  =&2\pi\sum_{\begin{subarray}{c}1\le j\le 2N\\y_j<y\end{subarray}}d_j=2\pi\sum_{j=1}^{2N}d_j\chi_{[y_j,1]}(y).
\end{align}
\eqref{eq:int j on [0,1]xx}, \eqref{eq:difference of int j} and \eqref{eq:degrees2} imply
\begin{align}\label{eq:intj1}
  \int_{\T}j_1(u(\vx))d\vx=&\int_0^1\int_0^1j_1(u(\vx))dxdy\nn
  =&\int_0^1\left(2\pi k_1-2\pi\sum_{j=1}^{2N}d_j\chi_{[y_j,1]}(y)\right)dy\nn
  =&2\pi k_1-2\pi\sum_{j=1}^{2N}d_j(1-y_j)=2\pi k_1+2\pi\sum_{j=1}^{2N}d_jy_j-2\pi\sum_{j=1}^{2N}d_j\nn
  =&\int_0^1j_1(u(x,0))dx+2\pi\sum_{j=1}^{2N}d_jy_j.
\end{align}
Similar to the proof of \eqref{eq:intj1}, there exists $k_2\in\Z$ such that
\begin{equation}\label{eq:intj2}
  \int_{\T}j_2(u(\vx))d\vx=2\pi k_2-2\pi\sum_{j=1}^{2N}d_jx_j=\int_0^1j_2(u(0,y))dy-2\pi\sum_{j=1}^{2N}d_jx_j.
\end{equation}
Combining \eqref{eq:intj1} and \eqref{eq:intj2}, we obtain \eqref{eq:limitation of momentum}. 

\end{proof}

Then, we will prove \eqref{eq:limit of momentum of u_02}, which is stated as the following lemma:
\begin{Lem}\label{lem:limitation of limit of current}
If $u^\varepsilon\in H^1(\T)$ is a sequence satisfying
\begin{equation}\label{eq:small limit of J}
  J(u^\varepsilon)\to\pi\sum_{j=1}^{2N}d_j\va_j \ \text{in}\ \M,,\quad
  E^\varepsilon(u^\varepsilon)\le 2N\pi\log\frac{1}{\varepsilon}+C,
\end{equation}
for some $\va=(\va_1,\cdots,\va_{2N})^T\in(\T)^{2N}_*$ and $C>0$, we have
\begin{equation}\label{eq:limitation of the limit of momentum}
  \lim_{\varepsilon\to 0}\int_{\T}\vj(u^\varepsilon)d\vx\in 2\pi\J\sum_{j=1}^{2N}d_j\va_j+2\pi\Z^2.
\end{equation}
\end{Lem}
\begin{proof}
Via \eqref{eq:small limit of J} and Theorem 1.4.3 in \cite{CollianderJerrard1999GLvorticesSE}, we have that $\|u^\varepsilon\|_{W^{1,3/2}(\T)}$ is uniformly bounded:
\begin{equation}
  \|u^\varepsilon\|_{W^{1,3/2}(\T)}\le C. 
\end{equation}
Hence, there exists a function $u\in W^{1,3/2}(\T)\hookrightarrow\hookrightarrow L^3(\T)$ such that
\begin{equation}\label{eq:conver of u in L3}
  u^\varepsilon\to u\ \text{in} \quad L^3(\T),\quad \nabla u^\varepsilon\wto \nabla u\ \text{in}\ L^{3/2}(\T),
\end{equation}
which together with \eqref{eq:def of e j J} implies 
\begin{equation}\label{eq:convergence of vj(ue) to vj(u)}
  \vj(u^\varepsilon)\wto \vj(u)\ \text{in}\ L^1(\T).
\end{equation}
Moreover, \eqref{eq:small limit of J} gives
\begin{equation}
  \left\||u^\varepsilon|^2-1\right\|_{L^2(\T)}^2\le 4\varepsilon^2 E^\varepsilon(u^\varepsilon)\to 0.
\end{equation}
Hence, we have $|u^\varepsilon|\to 1$, which together with \eqref{eq:conver of u in L3} implies that  $|u|=1$, i.e. $u$ is an $\mathbb{S}^1$-valued function. 

For any $\varphi\in C^1(\T)$, combining the definition of Jacobian \eqref{eq:def of e j J} and \eqref{eq:small limit of J}, we have
\begin{align*}
  \int_{\T}\varphi J(u)d\vx=&\frac{1}{2}\int_{\T}(\J\nabla\varphi)\cdot \vj(u)d\vx=\lim_{\varepsilon\to 0}\frac{1}{2}\int_{\T}(\J\nabla\varphi)\cdot \vj(u^\varepsilon)d\vx
  =\lim_{\varepsilon\to 0}\int_{\T}\varphi J(u^\varepsilon)d\vx=\pi\sum_{j=1}^{2N}d_j\varphi(\va_j),
\end{align*}
which implies $u$ is an $\mathbb{S}^1$-valued function satisfying \eqref{eq:def of Jacpobian}. Theorem \ref{eq:momentum of S1function} implies that $\int_{\T}\vj(u)d\vx$ satisfies \eqref{eq:limitation of momentum}. Since \eqref{eq:convergence of vj(ue) to vj(u)} holds, we have
\begin{equation}
  \int_{\T}\vj(u^\varepsilon)d\vx\to\int_{\T}\vj(u)d\vx,
\end{equation}
which together with \eqref{eq:limitation of momentum} implies \eqref{eq:limitation of the limit of momentum}.
\end{proof}

\subsection{Canonical harmonic map on torus}

 For any $\va=(\va_1,\cdots,\va_{2N})\in(\T)^{2N}_*$ and
$  \vq\in 2\pi\sum_{j=1}^{2N}d_j\va_j+2\pi\Z^2,$ Lemma 2.1 in \cite{ZhuBaoJian} and Lemma 2.2 in \cite{ZhuBaoJian} imply that there exists an $\mathbb{S}^1$-valued function $H=H(\vx)\in C^\infty_{loc}(\T\setminus\{\va_1,\cdots,\va_{2N}\};\mathbb{S}^1)\cap W^{1,1}(\T;\mathbb{S}^1)$ satisfies the following identities:
\begin{align}
\Div \vj(H)=0,\quad J(H)=\pi\sum_jd_j\delta_{\va_j},\quad \int_{\T}\vj(H)d\vx=&\J\bvec{q},\label{eq:div of jH}
\end{align}
and the function satisfying \eqref{eq:div of jH} is unique up to a phase transform.
Then we can define $H$ to be the canonical harmonic map with vortices $\va$. We will write $H=H(\vx;\va,\vq)$ to emphasize the dependence of $H$ on $\va$ and $\vq$.

\subsection{Renormalized energy on torus}

It should be noted that in the definition of the renormalized energy \eqref{eq:def of W}, $\vq$ satisfies $\vq\in2\pi\sum_{j=1}^{2N}d_j\va_j+2\pi\Z^2$. Hence, the derivative of $|\vq|^2$ with respect to $\va_j$ is not zero. We can always assume $\va_j=(x_j,y_j)^T\in(0,1)^2$ (If not, we can take a translation.). Then we can find a constant $\bvec{k}\in\Z^2$ such that locally, 
\begin{equation}
  \vq=2\pi\sum_{k=1}^{2N}d_k\va_k+2\pi\bvec{k},
\end{equation}
which together with \eqref{eq:def of W} implies
\begin{equation}\label{eq:expansion of W}
  W(\va;\vq)=-\pi\sum_{1\le k\ne l\le 2N}d_kd_lF(\va_k-\va_l)+\frac{1}{2}\left|2\pi\sum_{k=1}^{2N}d_k\va_j+2\pi\bvec{k} \right|^2.
\end{equation}
Taking the gradient with respect to $\va_j$ on both sides of \eqref{eq:expansion of W}, we have 
\begin{align}\label{eq:derivative of W}
  \nabla_{\va_j}W(\va;\vq)=&-2\pi\sum_{\begin{subarray}{c}1\le k\le 2N\\k\ne j\end{subarray}}d_kd_j\nabla F(\va_j-\va_k)+2\pi d_j\left(2\pi\sum_{k=1}^{2N}d_k\va_j+2\pi\bvec{k}\right)\nn
  =&2\pi d_j\left(-\sum_{\begin{subarray}{c}1\le k\le 2N\\k\ne j\end{subarray}}d_k\nabla F(\va_j-\va_k) +\vq\right).
\end{align}
Then for any small $\rho$, there exists a constant $C_\rho(\vq)$ such that $W$ is Lipschitz in $(\T)^{2N}_\rho$ where
\[
  (\T)^{2N}_\rho:=\{(\vx_1,\cdots,\vx_{2N})\in(\T)^{2N}||\vx_j-\vx_k|>\rho,\forall 1\le j<k\le 2N\},
\]
 i.e. there exists a constant $C_\rho$ such that:
\begin{equation}\label{eq:W is Lip}
 \| W(\va;\vq)\|_{C^1( (\T)^{2N}_\rho)}\le C_\rho.
\end{equation}

{For $\bvec{z}=(z_1,\cdots,z_K)^T,\bvec{w}=(w_1,\cdots,w_K)^T\in \mC^{K},K=1,2$, we define
\begin{equation}
  \la\bvec{z},\bvec{w}\ra=\re\sum_{j=1}^K\overline{z_j}w_j.
\end{equation}
In particular, if $\bvec{z},\bvec{w}\in \R^K\subset\mC^K$,
\begin{equation}
  \la\bvec{z},\bvec{w}\ra=\bvec{z}\cdot\bvec{w}.
\end{equation}
We denote $\Hess(\eta)$ the Hessian matrix of $\eta$.
And we define
\begin{equation}\label{eq:def of r}
  r(\va)=\frac{1}{4}\min_{1\le j<k\le 2N}|\va_j-\va_k|.
\end{equation}
For $\rho<r(\va)$, we define
\begin{equation}
  \T_{\rho}(\va):=\T\setminus\cup_{j=1}^{2N}B_\rho(\va_j), \quad \T_*(\va):=\cup_{\rho>0}\T_\rho(\va)=\T\setminus\{\va_1,\cdots,\va_{2N}\}.
\end{equation}}
Then Lemma 2.2 in \cite{ZhuBaoJian} implies 
\begin{equation}\label{eq:est of int of |jH|^2}
  \int_{\T_\rho(\va)}e^\varepsilon(\vj(H))d\vx=-2N\pi\log \rho+W(\va;\vq)+O(\rho^2),
\end{equation}
and Lemma 2.3 in \cite{ZhuBaoJian} implies that for any $\eta\in C^2(B_{r(\va)}(\va_j))$ which is linear in a neighborhood of $\va_{j}$, we have
\begin{align}
&\int_{\T} \la\Hess(\eta)\vj(H),\J\vj(H)\ra d \vx
=-\nabla \eta(\va_{j})\cdot(\J\nabla_{\va_{j}}W(\va;\vq)).\label{eq:production of jH and eta GP}
\end{align}

\subsection{Lemma related to dynamics on torus}
\begin{Lem}\label{lem:distosta}
There exists a universal constant $C$ such that for any sequence $u^\varepsilon\in H^1(\T\times[0,T_*])$ and Lipschitz functions $\va=\va(t)=(\va_1(t),\cdot,\va_{2N}(t))^T\in (\T)^{2N}_*$, $\vq(\va(t))\in 2\pi\sum_{j=1}^{2N}d_j\va_j(t)+2\pi\Z^2$, $u_*(\vx,t)=H(\vx;\va(t),\vq(\va(t)))$ satisfying
\begin{equation}\label{eq:contodirac1}
    J(u^\varepsilon(t))\to \pi\sum_{j=1}^{2N}d_j \delta _{\va_j(t)}\ \text{in}\ \M,\ \text{for}\ t\in[0,T_*],\end{equation}
    \begin{equation}\label{eq:converge of j/||to u_*}
     \frac{\vj(u^\varepsilon)}{|u^\varepsilon|}\wto\vj(u_*)\ \text{in}\ L^2_{loc}(\T_*(\va(t))\times[0,T_*]),
\end{equation}
and 
\begin{equation}\label{eq:ebcondition}
\Sigma(t):=\limsup_{\varepsilon\to 0}\left(\int_{\T}e^\varepsilon(u^\varepsilon(t))d\vx-W_\varepsilon(\va(t);\vq(\va(t)))\right)<+\infty,
\end{equation}
we have   for any $\rho<r_0:=\min_{t\in[0,T_*]}r(\va(t))$,
\begin{equation}\label{eq:distosta}
    \limsup_{\varepsilon\to 0}\int_0^{T_*}\int_{\T_{\rho}(\va)}\left(e^\varepsilon(|u^\varepsilon|)+\frac{1}{2}\left|\frac{\vj(u^\varepsilon)}{|u^\varepsilon|}-\vj(u_*) \right|^2\right)d\vx dt\le C \int_0^{T_*}\Sigma(t)dt.
\end{equation}
\end{Lem}

\begin{proof}
By Lemma 3 in \cite{JerrardSpirn2008RefinedEstimateGrossPitaevskiiVortices}, \eqref{eq:contodirac1} implies that for $\varepsilon<\rho$,
\begin{equation}\label{eq:lower bound of energy in balls}
  \int_{B_\rho(\va(t))}e^\varepsilon(u^\varepsilon(\vx,t))d\vx-\left(\gamma+\pi\log\frac{\rho}{\varepsilon}\right)\ge -\Sigma_1(\varepsilon),
\end{equation}
where 
\begin{equation}
  \Sigma_1(\varepsilon)=C\frac{\varepsilon}{\rho}\sqrt{\log\frac{\rho}{\varepsilon}}+\frac{C}{\rho}\sup_{t\in[0,T_*]}\|J(u^\varepsilon(t))-\pi\delta_{\vb_j(t)}\|_{\M}.
\end{equation}
\eqref{eq:lower bound of energy in balls}, \eqref{eq:ebcondition} and \eqref{eq:def of Wepsilon} imply
\begin{align}\label{eq:energy bound out balls}
  \int_{\T_{\rho}(\va(t))}e^\varepsilon(u^\varepsilon(\vx,t))d\vx=&\int_{\T}e^\varepsilon(u^\varepsilon(\vx,t))d\vx-\sum_{j=1}^{2N}\int_{B_\rho(\va_j(t))}e^\varepsilon(u^\varepsilon(\vx,t))d\vx\nn
\le&W_\varepsilon(\va(t);\vq(\va(t)))-2N\left(\gamma+\pi\log\frac{\rho}{\varepsilon}\right)+\Sigma(t)+\Sigma_1(\varepsilon)\nn 
  =&2N\pi\log\frac{1}{\rho}+W(\va(t);\vq(\va(t)))+\Sigma(t)+\Sigma_1(\varepsilon).
\end{align}
Combining \eqref{eq:est of int of |jH|^2} and \eqref{eq:energy bound out balls}, we get
\begin{equation}\label{eq:energy difference}
  \int_{\T_{\rho}(\va(t))}(e^\varepsilon(u^\varepsilon(t))-e^\varepsilon(u^*(t)))d\vx\le \Sigma(t)+\Sigma_1(\varepsilon)+O(\rho^2).
\end{equation}
\eqref{eq:def of e j J} implies
\begin{equation}\label{eq:differnce between eu eu*}
  e^\varepsilon(u^\varepsilon)-e^\varepsilon(u_*)=e^\varepsilon(|u^\varepsilon|)+\frac{1}{2}\left|\frac{\vj(u^\varepsilon)}{|u^\varepsilon|}-\vj(u_*) \right|^2+\vj(u_*)\cdot\left(\frac{\vj(u^\varepsilon)}{|u^\varepsilon|}-\vj(u_*)\right).
\end{equation}
Combining \eqref{eq:energy difference} and \eqref{eq:differnce between eu eu*}, we have
\begin{align}\label{eq:short time estimate}
  &\int_{\T_{\rho}(\va(t))}\left(e^\varepsilon(|u^\varepsilon(t)|)+\frac{1}{2}\left|\frac{\vj(u^\varepsilon(t))}{|u^\varepsilon(t)|}-\vj(u_*(t)) \right|^2\right)d\vx\nn&\quad\le \Sigma(t)+O(\rho^2)+\Sigma_1(\varepsilon)-\int_{\T_{\rho}(\va)}\vj(u_*(t))\cdot\left(\frac{\vj(u^\varepsilon(t))}{|u^\varepsilon(t)|}-\vj(u_*(t))\right)d\vx.
\end{align}

By \eqref{eq:converge of j/||to u_*}, integrating \eqref{eq:short time estimate} over $[0,T_*]$ on both sides and letting $\varepsilon\to 0$, we  obtain
\begin{equation}
  \limsup_{\varepsilon\to 0}\int_0^{T_*}\int_{\T_{\rho}(\va(t))}\left(e^\varepsilon(|u^\varepsilon|)+\frac{1}{2}\left|\frac{\vj(u^\varepsilon)}{|u^\varepsilon|}-\vj(u_*) \right|^2\right)d\vx dt \le \int_0^{T_*}\Sigma(t)dt+O(\rho^2).
\end{equation}

As the estimate above works for any $\rho'<\rho$ and $\T_{\rho',s}\supset\T_{\rho,s}$, we have
\begin{align}
  &\limsup_{\varepsilon\to 0}\int_0^{T_*}\int_{\T_{\rho}(\va)}\left(e^\varepsilon(|u^\varepsilon|)+\frac{1}{2}\left|\frac{\vj(u^\varepsilon)}{|u^\varepsilon|}-\vj(u_*) \right|^2\right)d\vx dt\nn&\quad \le\limsup_{\varepsilon\to 0}\int_0^{T_*}\int_{\T_{\rho'}(\va)}\left(e^\varepsilon(|u^\varepsilon|)+\frac{1}{2}\left|\frac{\vj(u^\varepsilon)}{|u^\varepsilon|}-\vj(u_*) \right|^2\right)d\vx dt\nn &\quad\le \int_0^{T_*}\Sigma(t)dt+O((\rho')^2).
\end{align}
Letting $\rho'\to0$, we obtain \eqref{eq:distosta}.
\end{proof}
\begin{Rmk}
Similar estimates were  obtained by Ignat-Jerrard in \cite{IgnatJerrard2021Renormalizedenergymanifold} (See Proposition 9.1 in \cite{IgnatJerrard2021Renormalizedenergymanifold}) for vector-valued functions $U:M\to TM$, where $M$ is a compact manifold and $TM$ is the tangent space of $M$. 
\end{Rmk}

\section[Converge of quantities related to ue]{Converge of quantities related to $u^\varepsilon$}\label{sec:CGL}
\subsection[Derivatives of quantities related to CGL]{Derivatives of quantities related to \texorpdfstring{\eqref{eq:CGL2}}{}}\label{sec:derivatives of quantities CGL}
Multiplying $\overline{u_t^\varepsilon}$ on both sides of \eqref{eq:CGL2} and taking the real part, we have
\begin{equation}\label{eq:derivative of e with respect to t CGL}
\frac{\p}{\p t}e^\varepsilon(u^\varepsilon(\vx,t))=\nabla\cdot\re(\overline{u^\varepsilon_t}\nabla u^\varepsilon)-k_\varepsilon|u^\varepsilon_t|^2,
\end{equation}
which implies that for any $0\le t_1<t_2$,
\begin{equation}
    E^\varepsilon(u^\varepsilon(t_2))-E^\varepsilon(u^\varepsilon(t_1))=-\int_{t_1}^{t_2} \int_{\T}k_\varepsilon|u_t^\varepsilon(\vx,\tau)|^2d\vx d\tau .\label{eq:energy decreasing CGL}
\end{equation}
Combining \eqref{eq:energy decreasing CGL} and \eqref{con:limit of integer of J(varphi)}, we obtain that there exists a $C>0$ such that for all $\varepsilon>0$,
\begin{equation}\label{con:weak bd of energy CGL}
E^\varepsilon(u^\varepsilon(t))\le 2N \pi\log\left(\frac{1}{\varepsilon} \right)+C.
 \end{equation}

Multiplying $\overline{u^\varepsilon}$ on both sides of \eqref{eq:CGL2} and taking the  imaginary part, we have
\begin{align}
\nabla\cdot \vj(u^\varepsilon(\vx,t))&=k_\varepsilon\im(\overline{u^\varepsilon}u_t^\varepsilon)+\frac{\lambda}{2}\frac{\p}{\p t}|u^\varepsilon|^2.\label{eq:divergence of j CGL}
\end{align}

{Multiplying $\nabla \overline{u^\varepsilon}$ on both sides of \eqref{eq:CGL2} and noting the fact \cite{CollianderJerrard1999GLvorticesSE}
\begin{equation}
    \p_t\vj(u^\varepsilon)=2\im (\overline{u^\varepsilon_t}\nabla u^\varepsilon)+\nabla\im(\overline{u^\varepsilon}u^\varepsilon_t),
\end{equation}
we have
\begin{equation}
    k_\varepsilon\re( \overline{u^\varepsilon_t}\nabla u^\varepsilon)+\frac{\lambda}{2}\frac{\p}{\p t}\vj(u^\varepsilon)=\nabla\cdot(\nabla u^\varepsilon\otimes\nabla u^\varepsilon)-\nabla e^\varepsilon(u^\varepsilon)+\frac{\lambda}{2}\nabla\im(\overline{u^\varepsilon}u^\varepsilon_t),\label{eq:derivative of j CGL}
\end{equation}}
which together with \eqref{eq:def of e j J} implies
that for any $\varphi\in C^\infty(\T)$,
\begin{equation}\label{eq:mainequalityCGL}
  -\int_{\T}k_\varepsilon\la\J\nabla u^\varepsilon,u_t^\varepsilon\nabla \varphi\ra d\vx+\lambda\int_{\T}\frac{\p}{\p t}J(u^\varepsilon)\varphi d\vx=\int_{\T}\la\Hess(\varphi)\nabla u^\varepsilon,\J\nabla u^\varepsilon\ra d\vx.
\end{equation}
In \eqref{eq:derivative of j CGL}, $\nabla u^\varepsilon\otimes \nabla u^\varepsilon$ and $\nabla\cdot (\nabla u^\varepsilon\otimes \nabla u^\varepsilon)$ follow the notation of (2.1.6) in \cite{CollianderJerrard1999GLvorticesSE}: $\nabla u^\varepsilon\otimes \nabla u^\varepsilon$ is the matrix
\begin{equation}
    \left(\begin{array}{cc}\la\p_xu^\varepsilon,\p_xu^\varepsilon\ra&\la\p_xu^\varepsilon,\p_yu^\varepsilon\ra\\\la\p_yu^\varepsilon,\p_xu^\varepsilon\ra&\la\p_yu^\varepsilon,\p_yu^\varepsilon\ra \end{array}\right)
\end{equation}
and 
\begin{equation}
    \nabla\cdot (\nabla u^\varepsilon\otimes \nabla u^\varepsilon)=\left(\begin{array}{c}\p_x\la\p_xu^\varepsilon,\p_xu^\varepsilon\ra+\p_y\la\p_yu^\varepsilon,\p_xu^\varepsilon\ra\\\p_x\la\p_xu^\varepsilon,\p_yu^\varepsilon\ra+\p_y\la\p_yu^\varepsilon,\p_yu^\varepsilon\ra \end{array}\right).
\end{equation}

\subsection[Existence and regularity of vortices: converge of Jacobian of ue]{Existence and regularity of vortices: converge of Jacobian of $u^\varepsilon$}
We first prove the existence and regularity of the vortices, which is summarized as the following lemma:
\begin{Lem}\label{thm: existence of vortices of cGL}
If the initial data $u^\varepsilon_0$ of \eqref{eq:CGL2} satisfies \eqref{con:limit of integer of J(varphi)} and \eqref{eq:mu=1}, then there exist Lipschitz paths $\vb_j:[0,T)\to\T$, such that
\begin{equation}\label{eq:convergence of J(u) CGL}
    J(u^\varepsilon(\vx,t))\to \pi\sum_{j=1}^{2N}d_j\delta_{\vb_j(t)}\ \text{in}\ \M,\quad k_\varepsilon e^\varepsilon(u^\varepsilon(\vx,t))\to\pi\sum_{j=1}^{2N}\delta_{\vb_j(t)}\ \text{in}\ \M.
\end{equation}
\end{Lem}
\begin{proof}
We denote $r(\va^0)$ by replacing $\va$ with $\va^0$ in \eqref{eq:def of r}.
For $\varepsilon$ small, assumption \eqref{con:limit of integer of J(varphi)} implies that
\begin{equation}\label{eq:distance between Ju_0 and dirad}
\left\|J(u_0^\varepsilon)- \pi\sum_{j=1}^{2N} d_j \delta _{\val_{j}} \right\|_{\M}\le \frac{\pi}{400}r(\va^0).
\end{equation}
For each small $\varepsilon$, we know $u^\varepsilon(t)\to u_0^\varepsilon$ as $t\to 0$ in $H^1(\T)$, so $J(u^\varepsilon(t))\to J (u_0^\varepsilon)$ in $L^1(\T)$ and also in $\M$. We define
\begin{equation}\label{eq:define of Te CGL}
T^\varepsilon:=\sup\left\{\tau>0:\|J(u_0^\varepsilon)- J(u^\varepsilon(t)) \|_{\M}\le \frac{\pi}{400}r(\va^0)\ \text{for any}\ t\in[0,\tau]\right\}.
\end{equation}
\eqref{eq:define of Te CGL} and \eqref{eq:distance between Ju_0 and dirad} imply that for $t<T^\varepsilon$, 
\begin{equation}\label{eq:distance between Ju and dirac}
\left\|J(u^\varepsilon(t))- \pi\sum_{j=1}^{2N} d_j \delta _{\val_{j}} \right\|_{\M}\le \frac{\pi}{200}r(\va^0)\ \text{for any}\ t<T^\varepsilon.
\end{equation}
\eqref{eq:distance between Ju and dirac}, \eqref{con:weak bd of energy CGL} and  Lemma 1.4.4 in \cite{CollianderJerrard1999GLvorticesSE} imply that there exist paths $\va_j^\varepsilon(t)\in\T, j=1,\cdots,2N,$ for each $t>0,\varepsilon<\varepsilon_0$ such that
\begin{align}
  &\left\|J(u^\varepsilon(\vx,t))-\pi\sum_{j=1}^{2N}d_j\delta_{\va_j^\varepsilon(t)} \right\|_{\M}=o(1),\quad\left\|k_\varepsilon e^\varepsilon(u^\varepsilon(\vx,t))-\pi\sum \delta_{\va_j^\varepsilon(t)}\right\|_{\M}=o(1),\label{eq:distance between Ju(t) and Dirace CGL}\\
     &E^\varepsilon(u^\varepsilon(t))\ge 2N\pi\log\frac{1}{\varepsilon}- C.\label{eq:lowerbound of E CGL}
\end{align}
Then \eqref{eq:lowerbound of E CGL} and \eqref{con:weak bd of energy CGL} imply that there exists a constant $C$ such that for any $t_1< t_2$,
\begin{equation}\label{eq:L2norm of u_t CGL}
  \int_{t_1}^{t_2} \int_{\T}k_\varepsilon|u_t^\varepsilon(\vx,\tau)|^2d\vx d\tau=E^\varepsilon(u^\varepsilon(t_2))-E^\varepsilon(u^\varepsilon(t_1))\le C.
\end{equation}

Then we  estimate $|\va_j^\varepsilon(t_2)-\va_j^\varepsilon(t_1)|$: we may find a function $\eta\in C_0^\infty(B_{r(\va^0)}(\va_j^0)) $ such that
\[
    \eta(\vx)=(\vx-\va_j^\varepsilon(t_1))\cdot\frac{\va^\varepsilon_j(t_2)-\va_j^\varepsilon(t_1)}{|\va^\varepsilon_j(t_2)-\va_j^\varepsilon(t_1)|},\quad \vx\in B_{3r(\va^0)/4}(\va_j^\varepsilon(t_1)).
\]
\eqref{eq:distance between Ju(t) and Dirace CGL}, \eqref{eq:L2norm of u_t CGL}, \eqref{con:weak bd of energy CGL} and \eqref{eq:def of E2} imply
\begin{align*}
&\pi|\va^\varepsilon_j(t_2)-\va_j^\varepsilon(t_1)|=\int_{\T}(\delta_{\va^\varepsilon_j(t_2)}-\delta_{\va^\varepsilon_j(t_1)})\eta d\vx\\
&\quad=\pi\int_{\T}k_\varepsilon(e^\varepsilon(u^\varepsilon(\vx,t_2))-e^\varepsilon(u^\varepsilon(\vx,t_1)))\eta(\vx) d\vx+o(1)\\
 &\quad\le k_\varepsilon\left|\int_{t_1}^{t_2}\int_{\T}\nabla\cdot \re(\overline{u^\varepsilon_t}\nabla u^\varepsilon)\eta d\vx dt \right|+k_\varepsilon^2\int_{t_1}^{t_2}\int_{\T}\eta|u_t^\varepsilon|^2d\vx dt +o(1)\\
&\quad\le k_\varepsilon\left\|\eta\right\|_{C^1(\T)}\|u_t^\varepsilon\|_{L^2(\T\times[t_1,t_2])}\|\nabla u^\varepsilon\|_{L^2(\T\times[t_1,t_2])}+o(1)\\
&\quad\le C\sqrt{t_2-t_1}+o(1).
\end{align*}
Hence, there is a lower bound $T_0$ of $\,T^\varepsilon$, and we can find continuous $\vb_j:[0,T_0]\to\T$ satisfying $\va_j^\varepsilon(t)\to \vb_j(t)$ up to a subsequence. As long as $\vb_j(t)\ne \vb_k(t)$ for any $j\ne k$, we can repeat the above proof and extend the existence interval of $\vb=(\vb_1,\cdots,\vb_{2N})^T$. The extension ends while reaching 
\begin{equation}\label{eq:def of T}
T=\min\{t:\va_j(t)=\va_k(t)\ \text{for some}\ j\ne k\}.
\end{equation}

It remains to prove that $\vb$ is Lipschitz. For any $T_1<T$, we denote $r=\min_{t\in[0,T_1]}r(\vb(t))$. For any $\rho<r$, via \eqref{eq:convergence of J(u) CGL}, the energy bound \eqref{con:weak bd of energy CGL} and Theorem 1.4.3 in \cite{CollianderJerrard1999GLvorticesSE}, we have
\begin{equation}
     \int_{\T_\rho(\vb(t))}e^\varepsilon(u^\varepsilon(\vx,t))d\vx \le C.\label{eq:upperboundofEout balls CGL}
 \end{equation} 
 For any $t_1,t_2$ such that $B_{r/2}(\vb_j(t))\subset B_{3r/4(\vb(t_1))}$ for any $t\in[t_1,t_2]$, we may find test functions $\varphi_1,\varphi_2\in C^\infty_0(B_{r(\va^0)}(\vb_j(t_1)))$ such that
\[
    \varphi_1(\vx)=(\vx-\vb_j(t_1))\cdot\bvec{\nu},\quad \varphi_2(\vx)=(\vx-\vb_j(t_1))\cdot\bvec{\nu}^\perp,\quad \vx\in B_{3r/4}(\vb_j(t_1)),
\]
where $\bvec{\nu}=d_j\lambda(\vb_j(t_2)-\vb_j(t_1))/|\lambda(\vb_j(t_2)-\vb_j(t_1))|$ and $\nu^\perp=-\J\nu$. 
Then  combining \eqref{eq:mainequalityCGL} and \eqref{eq:upperboundofEout balls CGL}, we have
\begin{align}\label{eq:estimate of |aj(t)-aj(s)|}
    &-k_\varepsilon\int_{t_1}^{t_2}\int_{\T}\la\J\nabla u^\varepsilon,u_t^\varepsilon\nabla \varphi_1\ra d\vx ds+\lambda\int_{t_1}^{t_2}\int_{\T}\varphi_1\frac{\p}{\p t}J(u^\varepsilon)d\vx ds\nn&\quad=\int_{t_1}^{t_2}\int_{\T}\la\Hess(\varphi_1)\nabla u^\varepsilon,\J\nabla u^\varepsilon\ra d\vx ds\nn
    &\quad\le|t_2-t_1|\|\varphi_1\|_{C^2(\T)}\sup_{t_1\le t\le t_2}\|\nabla u^\varepsilon\|^2_{L^2(\T_{r/2}(\vb(t)))}\le C|t_2-t_1|.
\end{align}
Via \eqref{eq:derivative of e with respect to t CGL}, \eqref{eq:L2norm of u_t CGL} and \eqref{eq:convergence of J(u) CGL}, and noting that $\nabla \varphi_2+\J\nabla\varphi_1=\bvec{0}$ in $B_{3r/4(\vb_j(t_1))}$, the first term in the first line of \eqref{eq:estimate of |aj(t)-aj(s)|} can be estimated as follows:
\begin{align}\label{eq:first term111}
&-k_\varepsilon\int_{t_1}^{t_2}\int_{\T}\la\J\nabla u^\varepsilon,u_t^\varepsilon\nabla \varphi_1\ra d\vx ds\nn
&\quad=-k_\varepsilon\int_{t_1}^{t_2}\int_{\T}\la\nabla u^\varepsilon,u^\varepsilon_t\nabla \varphi_2\ra d\vx ds+k_\varepsilon\int_{t_1}^{t_2}\int_{\T}\la\nabla u^\varepsilon,u^\varepsilon_t(\nabla \varphi_2+\J\nabla \varphi_1)\ra d\vx ds\nn
&\quad=k_\varepsilon\int_{t_1}^{t_2}\int_{\T}\varphi_2\frac{\p}{\p t}e^\varepsilon(u^\varepsilon)d\vx ds+k_\varepsilon^2\int_{t_1}^{t\bar{2}}\int_{\T}\varphi_2|u_t^\varepsilon|^2d\vx ds+O(\sqrt{k_\varepsilon})\nn
&\quad=(\vb_j(t+h)-\vb_j(t))\cdot\bvec{\nu}^\perp+o(1)=o(1).
\end{align}
The second term of the first line of \eqref{eq:estimate of |aj(t)-aj(s)|} converges to $\pi |\lambda(\vb_j(t+h)-\vb_j(t))|$ because of \eqref{eq:convergence of J(u) CGL}.
Hence, substituting \eqref{eq:first term111} into \eqref{eq:estimate of |aj(t)-aj(s)|} and letting $\varepsilon\to 0$, we have
\[
    \pi |\lambda(\vb_j(t+h)-\vb_j(t))|\le C|t_2-t_1|,
\]
which implies that $\vb$ is locally Lipschitz.

\end{proof}

\subsection[Converge of current of ue]{Converge of current of $u^\varepsilon$}
\begin{Lem}
Assume $u^\varepsilon$, $\vb$ and $T$ are the same as those in Lemma \ref{thm: existence of vortices of cGL}. Then for any $T_1<T$ and  $\rho<r$
\begin{equation}\label{eq:convergence of j(u) CGL}
    \vj(u^\varepsilon)\wto\vj(u^*)\ \text{in}\ L^1(\T\times[0,T_1]),\quad
    \frac{\vj(u^\varepsilon)}{|u^\varepsilon|}\wto\vj(u^*)\ \text{in}\ L^2(\T_\rho(\vb(t))\times[0,T_1]),
\end{equation}
 where 
 \begin{equation}\label{eq:def of u^*}u^*(\vx,t)=H(\vx;\vb(t),\vq(\vb(t))).\end{equation}
\end{Lem}

\begin{proof}
Similar to the proof of Lemma \ref{lem:limitation of limit of current}, we can find an $\mathbb{S}^1$-valued function $v\in L^\infty([0,T_1],W^{1,3/2}(\T))$ such that up to a subsequence,
\begin{equation}\label{eq:converddddd}
    \vj(u^\varepsilon)\wto \vj(v)\ \text{in}\ L^1(\T\times[0,T_1]),\quad J(v(t))=\pi\sum_{j=1}^{2N}d_j\delta_{\vb_j(t)}.
\end{equation}
Recall $r=\min_{t\in[0,T_1]}r(\vb(t))$. For $\rho<r$, \eqref{eq:upperboundofEout balls CGL} gives
 \begin{equation}
     \left\|\frac{\vj(u^\varepsilon)}{|u^\varepsilon|} \right\|^2_{L^2(\T_\rho(\vb(t))\times[0,T_1])}\le 2\int_0^{T_1}\int_{\T_\rho(\vb(t))}e^\varepsilon(u^\varepsilon(\vx,t))d\vx dt\le CT_1
 \end{equation}
 is uniformly bounded, which implies that there exists some $\vj_*\in L^2(\T_\rho(\vb(t))\times[0,T_1])$ such that up to a subsequence
 \begin{equation}\label{eq:raff estimate of j/||}
     \frac{\vj(u^\varepsilon)}{|u^\varepsilon|}\wto \vj_*\ \text{in}\ L^2(\T_\rho(\vb(t))\times[0,T_1]).
 \end{equation}
 Energy bound \eqref{con:weak bd of energy CGL} gives
 \begin{equation}\label{eq:bound of |u|-1}
 \left\||u^\varepsilon|-1\right\|_{L^2(\T)}\le C\varepsilon^2\log\frac{1}{\varepsilon},
 \end{equation}
 which implies $|u^\varepsilon|\to 1$ in $L^2(\T\times[0,T_1])$. Hence \eqref{eq:raff estimate of j/||} implies
\begin{equation}
    \vj(u^\varepsilon)=|u^\varepsilon|\frac{\vj(u^\varepsilon)}{|u^\varepsilon|}\wto \vj_*\ \text{in}\ L^1(\T_r(\va(t))\times[0,T_1]).
\end{equation}
Noting  \eqref{eq:converddddd}, we have $\vj_*(\vx,t)=\vj(v(\vx,t))$. So we only need to verify that $\vj(v(\vx,t))=\vj(u^*(\vx,t))$ to prove \eqref{eq:convergence of j(u) CGL}.

For any $\varphi\in C_0^\infty(\T\times[0,T_1])$, combining \eqref{eq:converddddd}, \eqref{eq:divergence of j CGL}, \eqref{eq:bound of |u|-1} and \eqref{eq:L2norm of u_t CGL}, we obtain
\begin{align}
    &\left|\int_{\T\times[0,T_1]}\nabla\varphi\cdot\vj(v)d\vx\right|=\left|\lim_{\varepsilon\to0}\int_{\T\times[0,T_1]}\nabla\varphi\cdot\vj(u^\varepsilon)d\vx\right|=\lim_{\varepsilon\to0}\left|\int_{\T\times[0,T_1]}\varphi\nabla\cdot\vj(u^\varepsilon)d\vx\right|\nn
    &\quad=\lim_{\varepsilon\to 0}\left|\int_{\T\times[0,T_1]} k_\varepsilon\varphi\im(\overline{u^\varepsilon}u_t^\varepsilon) d\vx +\int_{\T\times[0,T_1]}\varphi\frac{\p}{\p t}\frac{\lambda}{2}|u^\varepsilon|^2d\vx\right|\nn&\quad=\lim_{\varepsilon\to 0}\left|\int_{\T\times[0,T_1]} k_\varepsilon\varphi\im(\overline{u^\varepsilon}u_t^\varepsilon) d\vx -\int_{\T\times[0,T_1]}\varphi_t\frac{\lambda}{2}(|u^\varepsilon|^2-1)d\vx\right|\nn
    &\quad\le\lim_{\varepsilon\to 0}k_\varepsilon\|u_t^\varepsilon\|_{L^2(\T\times[0,T_1])} \|\varphi\|_{C(\T\times[0,T_1]}+\lim_{\varepsilon\to 0}\frac{|\lambda|}{2} \|\varphi\|_{C^1(\T\times[0,T_1])}\left\|u^\varepsilon|^2-1\right\|_{L^2(\T\times[0,T_1])} =0,
\end{align}
which implies $\nabla\cdot \vj(v)=0$.
Noting \eqref{eq:def of u^*} and \eqref{eq:div of jH}, we have
\[
    \nabla \cdot (\vj(v)-\vj(u^*))=0.
\]
Similarly, we have
\[
  \nabla\cdot(\J(\vj(v)-\vj(u^*)))=0.
\]
Hence, there exists a function $\bvec{g}:[0,T_1]\to \R^2$ such that $\vj(v(\vx,t))-\vj(u^*(\vx,t))=\bvec{g}(t)$.  We have  $\bvec{g}(0)=\bvec{0}$ via the definition of $u^*$ \eqref{eq:def of u^*}, the definition of $v$ \eqref{eq:converddddd}, \eqref{con:limit of integer of J(varphi)} and \eqref{eq:div of jH}. Then, since both $v$ and $u^*$ are $\mathbb{S}^1$-valued functions on torus, Lemma \ref{eq:momentum of S1function} together with  \eqref{eq:def of u^*}, \eqref{eq:div of jH} and \eqref{eq:converddddd} implies
\[
    \bvec{g}(t)=\int_{\T}\vj(v(\vx,t))d\vx-\int_{\T}\vj(u^*(\vx,t))d \vx\in 2\pi\Z^2.
\]
So, we only need to verify the continuity of $\bvec{Q}(v(t))=\int \vj(v(\vx,t))d\vx $ to derive that $\bvec{g}=\bvec{0}$.
\eqref{eq:derivative of j CGL} implies 
\begin{align*}
    \left|\frac{\lambda}{2}(\bvec{Q}(v(t_2))-\bvec{Q}(v(t_1)))\right|=&\lim_{\varepsilon\to 0}\frac{|\lambda|}{2}\left|\int_{\T}(\vj(u^\varepsilon(\vx,t_2))-\vj(u^\varepsilon(\vx,t_1)))d\vx \right|\\
    =&\lim_{\varepsilon\to 0}\left|\int_{t_1}^{t_2}dt\int_{\T}k_\varepsilon u_t^\varepsilon\cdot \nabla u^\varepsilon d\vx \right|\\
    \le&C\limsup_{\varepsilon\to 0}k_\varepsilon\|u_t^\varepsilon\|_{L^2(\T\times[t_1,t_2])}\|\nabla u^\varepsilon\|_{L^2(\T\times[t_1,t_2])}\\\le& C\sqrt{t_2-t_1}.
\end{align*} 
Hence, $\bvec{Q}(v(t))$ is continuous, which implies \eqref{eq:convergence of j(u) CGL} immediately. 
\end{proof}

\section{Dynamics of vortices: proof of Theorem \ref{thm:dynamic CGL2}}\label{sec:dynamics}
\begin{proof}
Assume $\va=\va(t)=(\va_1(t),\cdots,\va_{2N}(t))^T$ is the solution to equation \eqref{eq:CGLODE} and recall $\vb=\vb(t)=(\vb_1(t),\cdots,\vb_{2N}(t))^T$ is obtained in Lemma \ref{thm: existence of vortices of cGL}. We define
\begin{equation}\label{eq:def of zeta}
    \zeta_j(t):=\sqrt{1+\lambda^2}|\va_j(t)-\vb_j(t) |=|(I_2-\lambda d_j\J)(\va_j(t)-\vb_j(t))|,
\end{equation}
where $I_2$ is the unit matrix of second order
\[
  I_2=\left(\begin{array}{cc}1&0\\0&1 \end{array}\right).
\]
For simplicity, we will use the notation in this subsection: 
\begin{equation}\label{eq:simple notation}
\begin{array}{ll}
    W(\va(t))=W(\va(t),\vq(\va(t))),&\quad W(\vb(t))=W(\vb(t),\vq(\vb(t))),\\
    W_\varepsilon(\va(t))=W_\varepsilon(\va(t),\vq(\va(t))),&\quad W_\varepsilon(\vb(t))=W_\varepsilon(\vb(t),\vq(\vb(t))).
\end{array}
\end{equation}
We can find $T_2<T$ such that for any $t<T_2$, we have
\begin{equation}\label{eq:def of r2 and T2}
  \max_{1\le j\le 2N}\zeta_j(t)<r_2:=\inf_{t\in[0,T_2]}\min\{r(\va(t)),r(\vb(t))\}.
\end{equation}
For any $t\in[0,T_2]$, taking the derivative of $\zeta$ with respect to $t$, we obtain
\begin{align}\label{eq:estimate of dzeta}
|\dot{\zeta}_j(t)|\le&|(I_2-\lambda d_j\J)(\dot{\va}_j(t)-\dot{\vb}_j(t))|\nn
\le&\frac{1}{\pi}\left|\nabla_{\va_j}W(\va(t))-\nabla_{\vb_j}W(\vb(t)) \right|+\left|(I_2-\lambda d_j\J)\dot{\vb}_j(t)+\frac{1}{\pi} \nabla_{\vb_j}W(\vb(t)) \right|\nn
\le& C\sum_{j=1}^{2N}\zeta_j(t)+\left|(I_2-\lambda d_j\J)\dot{\vb}_j(t)+\frac{1}{\pi} \nabla_{\vb_j}W(\vb(t)) \right|.
\end{align}
The last inequality holds because of \eqref{eq:def of r2 and T2} and the Lipschitz property of $W$ \eqref{eq:W is Lip}.

We take $\bvec{\nu}$ satisfying
\begin{equation}\label{eq:def of nu}
    \bvec{\nu}^\perp:=-\J\bvec{\nu}=\frac{(I_2-\lambda d_j\J)\dot{\vb}_j(t)+\frac{1}{\pi} \nabla_{\vb_j}W(\vb(t))}{ \left|(I_2-\lambda d_j\J)\dot{\vb}_j(t)+\frac{1}{\pi} \nabla_{\vb_j}W(\vb(t))\right|},
\end{equation}
and find $\varphi_1,\varphi_2\in C^\infty(B_{r_2}(\vb_j(t)))$ satisfying 
\begin{equation}\label{eq:def of phi}
    \varphi_1(\vx)=(\vx-\vb_j(t))\cdot\bvec{\nu},\quad \varphi_2(\vx)=(\vx-\vb_j(t))\cdot\bvec{\nu}^\perp,\quad \vx\in B_{3r_3/4}(\vb_j(t)).
\end{equation}
\eqref{eq:def of phi} immediately implies
\begin{equation}\label{eq:Df2+JDf1}
  \nabla \varphi_2+\J\nabla\varphi_1=0\ \text{in}\ B_{3r_3/4}(\vb_j(t)).
\end{equation}
 Substituting $\varphi_1$ into \eqref{eq:mainequalityCGL}, we have
\begin{equation}\label{eq:CGLmiddle}
    -\int_{\T}k_\varepsilon\la\J\nabla u^\varepsilon,u_t^\varepsilon\nabla \varphi_1\ra d\vx+\lambda\int_{\T}\frac{\p}{\p t}J(u^\varepsilon)\varphi_1 d\vx=\int_{\T}\la\Hess(\varphi_1)\nabla u^\varepsilon,\J\nabla u^\varepsilon\ra d\vx.
\end{equation}
Via \eqref{eq:derivative of e with respect to t CGL}, the first term on the left hand side of \eqref{eq:CGLmiddle} can be estimated as follows:
\begin{align}\label{eq:fisrt term of 9}
    &-\int_{\T}k_\varepsilon\la\J\nabla u^\varepsilon,u_t^\varepsilon\nabla \varphi_1\ra d\vx\nn
    &\quad=-\int_{\T}k_\varepsilon\la\nabla u^\varepsilon,u_t^\varepsilon\nabla \varphi_2\ra d\vx+k_\varepsilon\int_{\T}\la\nabla u^\varepsilon,u_t^\varepsilon\nabla \varphi_2\ra+\J\nabla\varphi_1)d\vx\nn
    &\quad=k_\varepsilon\int_{\T}\varphi_2\frac{\p}{\p t}e^\varepsilon(u^\varepsilon)d\vx+k_\varepsilon^2\int_{\T}\varphi_2|u^\varepsilon_t|^2d\vx+k_\varepsilon\int_{\T}\la\nabla u^\varepsilon,u_t^\varepsilon(\nabla \varphi_2 +\J\nabla\varphi_1)\ra d\vx.
\end{align}
Substituting \eqref{eq:fisrt term of 9} into \eqref{eq:CGLmiddle}, we have
\begin{align}\label{eq:estimate of CGL 2}
    &k_\varepsilon\int_{\T}\varphi_2\frac{\p}{\p t}e^\varepsilon(u^\varepsilon)d\vx+\lambda\int_{\T}\frac{\p}{\p t}J(u^\varepsilon)\varphi_1 d\vx+k_\varepsilon^2\int_{\T}\varphi_2|u^\varepsilon_t|^2d\vx+k_\varepsilon\int_{\T}\la\nabla u^\varepsilon,u_t^\varepsilon(\nabla \varphi_2+\J\nabla\varphi_1)\ra d\vx\nn 
    &\quad=\int_{\T}\la\Hess(\varphi_1)\nabla u^\varepsilon,\J\nabla u^\varepsilon\ra d\vx.
\end{align}
Integrating \eqref{eq:estimate of CGL 2} on $[t,t+h]$, then letting $\varepsilon\to 0$ and using \eqref{eq:convergence of J(u) CGL},  \eqref{eq:upperboundofEout balls CGL}, \eqref{eq:L2norm of u_t CGL} and \eqref{eq:Df2+JDf1}, we finally have 
\begin{align}\label{eq:estimate of b(t+h)-b(t)}
    &\bvec{\nu}^\perp\cdot(\vb_j(t+h)-\vb_j(t))+d_j\lambda\bvec{\nu}\cdot(\vb_j(t+h)-\vb_j(t))=\bvec{\nu}^\perp\cdot((I_2-\lambda d_j\J)(\vb_j(t+h)-\vb_j(t)))\nn&\quad=\lim_{\varepsilon\to 0}\int_t^{t+h}\int_{\T}\la\Hess(\varphi_1)\nabla u^\varepsilon,\J\nabla u^\varepsilon\ra d\vx ds.
\end{align}
Combining \eqref{eq:def of u^*}, \eqref{eq:production of jH and eta GP}, \eqref{eq:def of nu}, \eqref{eq:def of phi} and \eqref{eq:estimate of b(t+h)-b(t)}, we get
\begin{align}\label{eq:splitting of estimate}
&\left|(I_2-\lambda d_j\J)\dot{\vb}(t)+\frac{1}{\pi} \nabla_{\vb_j}W(\vb(t)) \right|=\bvec{\nu}^\perp\cdot\left(\lim_{h\to 0}(I_2-\lambda d_j\J)\frac{\vb_j(t+h)-\vb_j(t)}{h}+\frac{1}{\pi}\nabla_{\vb_j}W(\vb(t)) \right)\nn
&\quad=\lim_{h\to0}\lim_{\varepsilon\to0}\frac{1}{\pi h}\int_{t}^{t+h}\int_{\T}\left(\la\Hess(\varphi_1)\nabla u^\varepsilon,\J\nabla u^\varepsilon\ra-\la\Hess(\varphi_1)\vj(u^*),\J\vj(u^*)\ra \right)d\vx ds\nn
&\quad=L_j+K_{1j}+K_{2j}+K_{3j} ,
\end{align}
with
\begin{align}
    L_j(t)=&\lim_{h\to 0}\lim_{\varepsilon\to 0}\frac{1}{\pi h}\int_t^{t+h}\int_{\T}\la\Hess(\eta)\nabla|u^{\varepsilon}|,\J\nabla|u^{\varepsilon}| \ra d \vx\ud s,\label{eq:define of L}\\
K_{j1}(t)=&\lim_{h\to 0}\lim_{\varepsilon\to 0}\frac{1}{\pi h}\int_t^{t+h}\int_{\T}\la\Hess(\eta)\left(\frac{\vj(u^{\varepsilon})}{|u^{\varepsilon}|}-\vj(u_*)\right),\J\left(\frac{\vj(u^{\varepsilon})}{|u^{\varepsilon}|}-\vj(u_*)\right)\ra d \vx\ud s,\label{eq: 1st term of L} \\
K_{j2}(t)=&\lim_{h\to 0}\lim_{\varepsilon\to 0}\frac{1}{\pi h}\int_t^{t+h}\int_{\T}\la\Hess(\eta)\vj(u_*),\J\left(\frac{\vj(u^{\varepsilon})}{|u^{\varepsilon}|}-\vj(u_*)\right)\ra d \vx\ud s,\label{eq: 2nd term of L} \\
K_{j3}(t)=&\lim_{h\to 0}\lim_{\varepsilon\to 0}\frac{1}{\pi h}\int_t^{t+h}\int_{\T}\la\Hess(\eta)\left(\frac{\vj(u^{\varepsilon})}
{|u^{\varepsilon}|}-\vj(u_*)\right),\J\vj(u_*)\ra d \vx\ud s,\label{eq: 3rd term of L}
{}\end{align}
where we adopt
\begin{equation}\label{eq:separation of nabla square}
  \begin{cases}
  \la\p_xu^\varepsilon,\p_xu^\varepsilon\ra&=\frac{1}{|u^\varepsilon|^2}(j_1(u^\varepsilon))^2+(\p_x|u^\varepsilon|)^2,\\
  \la\p_xu^\varepsilon,\p_yu^\varepsilon\ra&=\frac{1}{|u^\varepsilon|^2}j_1(u^\varepsilon)j_2(u^\varepsilon)+\p_x|u^\varepsilon|\,\p_y|u^\varepsilon|,\\
  \la\p_yu^\varepsilon,\p_yu^\varepsilon\ra&=\frac{1}{|u^\varepsilon|^2}(j_2(u^\varepsilon))^2+(\p_y|u^\varepsilon|)^2, 
  \end{cases}
\end{equation}
which is a corollary of 
\begin{equation}
  \nabla u^\varepsilon=\frac{u^\varepsilon}{|u^\varepsilon|}\nabla|u^\varepsilon|+\frac{\vj(u^\varepsilon)}{|u^\varepsilon|}\frac{iu^\varepsilon}{|u^\varepsilon|}.
\end{equation}
Then \eqref{eq:convergence of j(u) CGL} and \eqref{eq:def of phi} imply that
\begin{equation}\label{eq:Kj2+Kj3}
K_{j2}=0,K_{j3}=0.
\end{equation}

Corollary 7 in \cite{SandierSerfaty2004EstimateGinzburgLandau}, \eqref{eq:energy decreasing CGL}  and \eqref{eq:convergence of J(u) CGL} imply
\begin{equation}\label{eq:uppper bound of Et-E0}
    E^\varepsilon(u^\varepsilon(t))=E^\varepsilon(u^\varepsilon_0)-\int_0^tds\int_{\T}k_\varepsilon|\p_tu^\varepsilon|^2d\vx \le E^\varepsilon(u^\varepsilon_0)-\pi\int_0^t|\dot{\vb}(s)|^2ds.
\end{equation}
Substituting \eqref{eq:derivative of W} and \eqref{eq:CGLODE} into the derivative of $W(\va(t))$ with respect to $t$, we obtain
\begin{align*}
\frac{d}{dt}W(\va(t))=&\sum_{j=1}^{2N}\nabla_{\va_j}W(\va(t))\cdot\dot{\va}_j(t)=-\pi\sum_{j=1}^{2N}(\dot{\va}_j(t)-\lambda d_j\J\dot{\va}_j(t))\cdot\dot{\va}_j(t)=-\pi |\dot{\va}_j(t)|^2,
\end{align*}
which implies
\begin{equation}\label{eq:conser of W(va)}
\pi\int_0^t|\dot{\va}(s)|^2ds+W(\va(t))\equiv W(\va(0)).
\end{equation}
Combining \eqref{eq:uppper bound of Et-E0}, \eqref{con:limit of integer of J(varphi)}, \eqref{eq:conser of W(va)}, \eqref{eq:def of Wepsilon} and \eqref{eq:W is Lip}, we have
\begin{align}\label{eq:strong energy bound of Et}
&\limsup_{\varepsilon\to 0}\left(E^\varepsilon(u^\varepsilon(t))-W_\varepsilon(\vb(t))\right)\le \lim_{\varepsilon\to 0}\left( W_\varepsilon(\va^0)-\pi\int_0^t|\dot{\vb}(s)|^2ds-W_\varepsilon(\vb(t))\right)\nn
&\quad=W(\va(t))+\int_0^t\left(|\dot{\va}(s)|^2-|\dot{\vb}(s)|^2\right)ds-W(\vb(t))\nn
&\quad\le C\sum_{j=1}^{2N}\zeta_j(t)+\int_0^t\sum_{j=1}^{2N}|\dot{\zeta}_j(s)|ds.
\end{align}
Lemma \ref{lem:distosta} and \eqref{eq:strong energy bound of Et} implies that 
\begin{equation}\label{eq:estimate e(||)+|j/||-j|}
    \int_{\T_{3r_2/4}(\vb(s))}\left(e^\varepsilon(|u^\varepsilon|)+\left|\frac{\vj(u^\varepsilon)}{|u^\varepsilon|}-\vj(u^*) \right|^2\right)d\vx\le C\sum_{j=1}^{2N}\zeta_j(t)+\int_0^t\sum_{j=1}^{2N}|\dot{\zeta}_j(s)|ds.
\end{equation}
 Substituting \eqref{eq:def of phi}, \eqref{eq:estimate e(||)+|j/||-j|} and \eqref{eq:Kj2+Kj3} into \eqref{eq:splitting of estimate}, we have
\begin{equation}\label{eq:estimate of second part}
    \left|(I_2-\lambda d_j\J)\dot{\vb}(t)+\frac{1}{\pi} \nabla_{\vb_j}W(\vb(t)) \right|=K_{1j}+L_j\le C\sum_{j=1}^{2N}\zeta_j(t)+\int_0^t\sum_{j=1}^{2N}|\dot{\zeta}_j(s)|ds.
\end{equation}
Substituting \eqref{eq:estimate of second part} into \eqref{eq:estimate of dzeta}, we have
\begin{equation}\label{eq:estimate of z' final}
    \sum_{j=1}^{2N}|\dot{\zeta}_j(t)|\le C\sum_{j=1}^{2N}\zeta_j(t)+\int_0^t\sum_{j=1}^{2N}|\dot{\zeta}_j(s)|ds\ \text{in}\ [0,T_2].
\end{equation}
Since $\zeta_j(0)=0$ which is a corollary of \eqref{eq:initial data CGL}, \eqref{eq:convergence of J(u) CGL}, \eqref{con:limit of integer of J(varphi)} and \eqref{eq:def of zeta}, \eqref{eq:estimate of z' final} implies
\[
    \zeta_j(t)=0\ \text{in}\ [0,T_2],
\]
which means that $\vb=\va$ in $[0,T_2]$. In particular, we have $\zeta_j(T_2)=0$ for any $1\le j\le 2N$. Then we can repeat the proof above to extend the interval such that $\zeta_j(t)=0$ holds to $[0,T)$, i.e. $\va(t)\equiv \vb(t)$ on $[0,T)$. Combining the definition of $\va$ \eqref{eq:CGLODE} and \eqref{eq:convergence of J(u) CGL}, we have proved Theorem \ref{thm:dynamic CGL2}.
\end{proof}

\section{Properties of \eqref{eq:CGLODE}}\label{sec:properties of CGLODE}
\subsection{First integrals}
Define
\begin{equation}\label{eq:xi(va)}
  \bvec{\xi}(\va):=\sum_{j=1}^{2N}\va_j-\lambda\J\sum_{j=1}^{2N}d_j\va_j,\quad \va\in(\T)^{2N}_*.
\end{equation}
\begin{Thm}
Assume $\va=(\va_1\cdots, \va_{2N})^T$ is the solution of \eqref{eq:CGLODE} with initial data \eqref{eq:initial data CGL}, then $\bvec{\xi}(\va)$ defined in \eqref{eq:xi(va)} is a first integral of \eqref{eq:CGLODE}, i.e. 
\begin{equation}\label{eq:firs integral}
  \bvec{\xi}(\va):=\bvec{\xi}(\va(t))\equiv\bvec{\xi}(\va(0))=\bvec{\xi}(\va^0).
\end{equation}
\end{Thm}
\begin{proof}
Differentiating \eqref{eq:xi(va)} with respect to $t$, and noting \eqref{eq:CGLODE}, \eqref{eq:degrees2} and \eqref{eq:derivative of W}, and that $F$ is an even function, we have
\begin{align}
  \frac{d}{dt}\bvec{\xi}(\va(t))&=\sum_{j=1}^{2N}\left(\dot{\va}_j(t)-\lambda d_j\J\dot{\va}_j(t) \right)=-\sum_{j=1}^{2N}\frac{1}{\pi}\nabla_{\va_j}W(\va(t);\vq(\va(t)))\nn 
  &=\sum_{j=1}^{2N}2d_j\left(\sum_{\begin{subarray}{c}1\le k\le 2N\\k\ne j \end{subarray}}d_k\nabla F(\va_j(t)-\va_k(t))-\vq(\va(t)) \right)\nn 
  &=2\sum_{1\le j\ne k\le 2N}d_jd_k\nabla F(\va_j(t)-\va_k(t))-\left(\sum_{j=1}^{2N}d_j\right)\vq(\va(t))=\bvec{0},
\end{align}
which immediately implies \eqref{eq:firs integral}.
\end{proof}

\subsection{Solutions for several initial data with symmetry}\label{sec:solutions}

\begin{Lem}
If $N=1$ and initial data \eqref{eq:initial data CGL} satisfies
\begin{equation}\label{eq:intial data 2v}
  \va^0_1=(0.5,0.5)^T+(\alpha^0,\beta^0)^T,\quad \va_2^0=(0.5,0.5)^T+(-\alpha^0,\beta^0),
\end{equation}
 then the solution $\va=\va(t)=(\va_1(t),\va_2(t))$ is of the following form
\begin{equation}\label{eq:symmetry of a}
  \va_1(t)=(0.5,0.5)^T+(\alpha(t),\beta(t))^T,\quad \va_2(t)=(0.5,0.5)^T+(-\alpha(t),\beta(t))^T, 
\end{equation}
where $\alpha,\beta$ satisfy
\begin{equation}\label{eq:sol of alpha beta}
  \dot{\alpha}=\frac{-2(\p_xF(2\alpha,0)+4\pi\alpha)}{1+\lambda^2},\quad\dot{\beta}=\frac{2\lambda(\p_xF(2\alpha,0)+4\pi\alpha)}{1+\lambda^2},
\end{equation}
with initial data
\begin{equation}
  \alpha(0)=\alpha^0,\beta(0)=\beta^0.
\end{equation}
In particular, the trajectory of solution $\va$ consist of two segments.

\end{Lem}
\begin{proof}
By the symmetry of \eqref{eq:define of F}, we have $F$ satisfies
\begin{equation}\label{eq:symmetry of F}
  F(x,y)=F(-x,y)=F(x,-y)=F(y,x).
\end{equation}
Noting that $F\in C^\infty_0(\T\setminus\{\bvec{0}\})$, we have
\begin{equation}\label{eq:derivative of F on edge}
  \p_x(0.5,y)=0,\quad \p_y(x,0.5)=0, \quad \text{for}\  x,y\in[0,1],
\end{equation}
and
\begin{equation}\label{eq:derivative of F on middle}
  \p_x(0,y)=\p_x(1,y)=0,\quad \p_y(x,0)=\p_y(x,1)=0, \quad \text{for}\  x,y\in  (0,1).
\end{equation}
Then the symmetry of initial data \eqref{eq:intial data 2v} and the symmetry of \eqref{eq:CGLODE}, we can take the ansatz that the solution $\va$  satisfies \eqref{eq:symmetry of a}. Substituting \eqref{eq:symmetry of a} into \eqref{eq:CGLODE} and noting \eqref{eq:derivative of F on middle},  \eqref{eq:derivative of W}, we have
\begin{equation}
  \left(\begin{array}{cc}1&-\lambda\\\lambda&1\end{array}\right)\dot{\va}_1=-2\nabla F(\va_1-\va_2)-4\pi(\va_1-\va_2)=2\left(\begin{array}{c}-\p_xF(2\alpha,0)-4\pi\alpha\\0\end{array}\right),
\end{equation}
which together with $\dot{\va}_1=(\dot{\alpha},\dot{\beta})^T$ immediately implies \eqref{eq:sol of alpha beta}.
Since $\dot{\alpha}/\dot{\beta}\equiv -1/\lambda$, \eqref{eq:symmetry of a} implies that the trajectory of $\va$ consists of two segments with slopes $-1/\lambda$ and $1/\lambda$.
\end{proof}
\begin{figure}[htp!]
\begin{center}
\begin{tabular}{cccc}
\includegraphics[height=6cm]{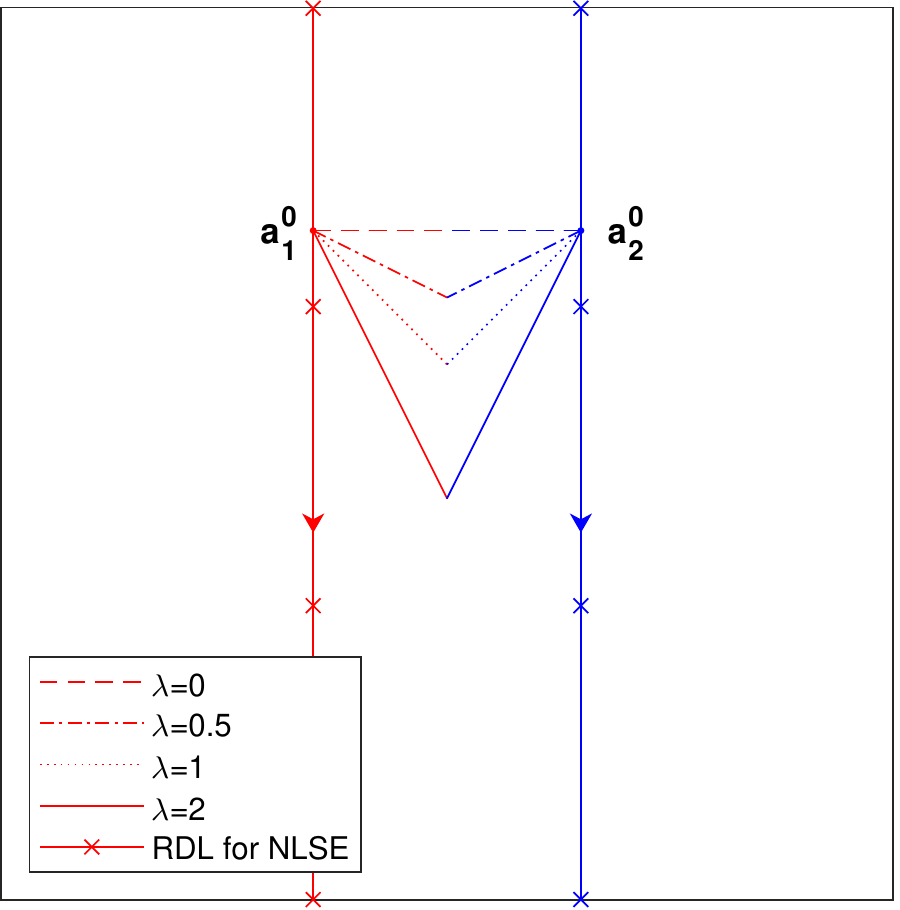}
&\includegraphics[height=6cm]{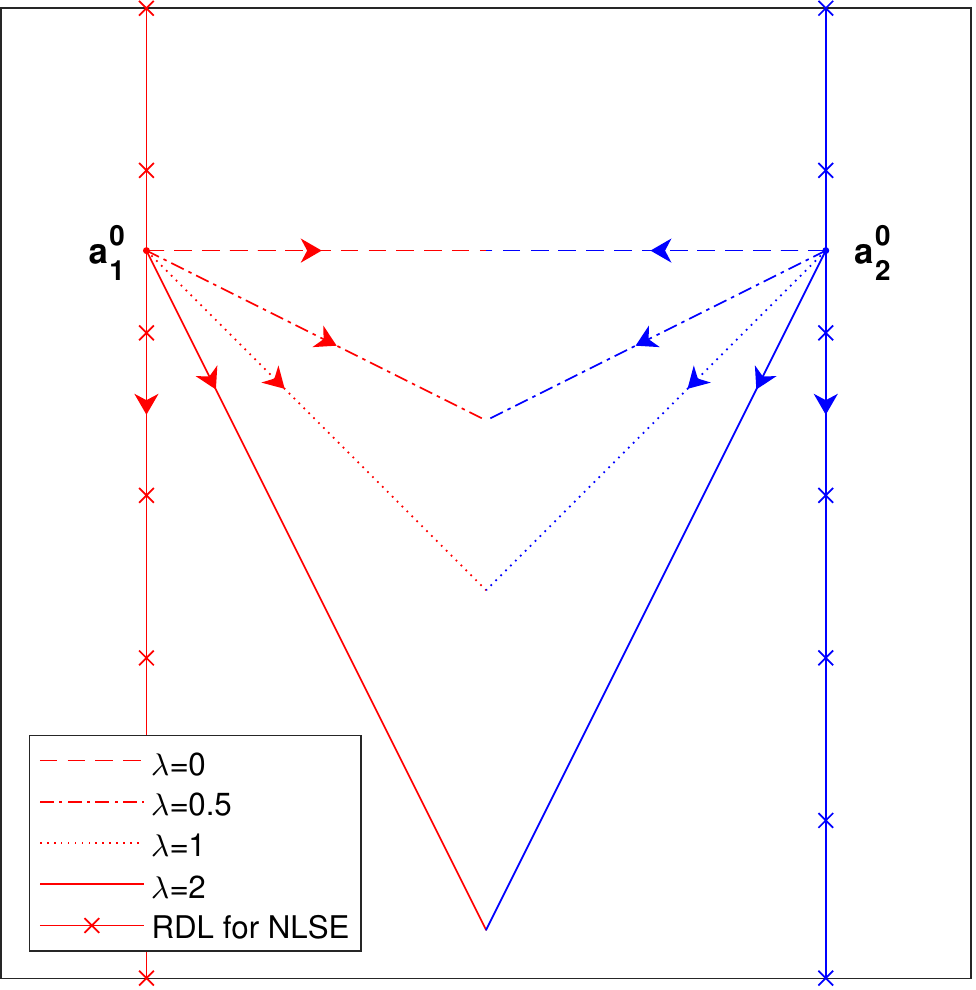}\\
$(\alpha^0,\beta^0)=(-0.15,0.25)$&$(\alpha^0,\beta^0)=(-0.35,0.25)$
\end{tabular}
\end{center}
\caption{Some trajectories of some solution of \eqref{eq:CGLODE} with different $\lambda$, the reduced dynamical laws of \eqref{eq:NLS} (marked by ``RDL for NLSE'') and the reduced dynamical laws of \eqref{eq:GL2} (marked by ``$\lambda=0$'') with initial data of the type \eqref{eq:intial data 2v}. 
}\label{fig:path 2 v}
\end{figure}
Figure \ref{fig:path 2 v} shows some numerical results for the solution of \eqref{eq:CGLODE} with initial data satisfying \eqref{eq:symmetry of a}. It can be seen from Figure \ref{fig:path 2 v} that  the trajectory of the solution of \eqref{eq:CGLODE} converges to  the trajectory of the solution of the reduced dynamical law of \eqref{eq:GL2} as $\lambda\to 0$, and converges to the trajectory of the solution of the reduced dynamical law of \eqref{eq:NLS} \cite{CollianderJerrard1999GLvorticesSE,ZhuBaoJian} as $\lambda\to \infty$.

\begin{Lem}
If $N=2$ and initial data \eqref{eq:initial data CGL} satisfies
\begin{equation}\label{eq:intial data 4v}
\begin{aligned}
  \va^0_1=(0.5,0.5)^T+(\alpha^0,\beta^0)^T,\quad\va_2^0=(0.5,0.5)^T+(-\alpha^0,-\beta^0),\\ \va^0_3=(0.5,0.5)^T+(-\alpha^0,\beta^0)^T,\quad\va_4^0=(0.5,0.5)^T+(\alpha^0,-\beta^0),
  \end{aligned}
\end{equation}
 then the solution $\va=\va(t)=(\va_1(t),\va_2(t),\va_3(t),\va_4(t))$ is of the following form
\begin{equation}\label{eq:symmetry of a 4v}
\begin{aligned}
  \va_1(t)=(0.5,0.5)^T+(\alpha(t),\beta(t))^T,\quad\va_2(t)=(0.5,0.5)^T+(-\alpha(t),-\beta(t)),\\\va_3(t)=(0.5,0.5)^T+(-\alpha(t),\beta(t))^T,\quad\va_4(t)=(0.5,0.5)^T+(\alpha(t),-\beta(t)),
  \end{aligned}
\end{equation}
where $\alpha,\beta$ satisfy
\begin{equation}\label{eq:sol of alpha beta 4v}
  \left(\begin{array}{c}\dot{\alpha}\\\dot{\beta} \end{array}\right)=\frac{2}{1+\lambda^2}\left(\begin{array}{cc}1&\lambda\\-\lambda&1\end{array}\right)\left(\begin{array}{c}\p_xF(2\alpha,2\beta)-\p_xF(2\alpha,0)\\\p_yF(2\alpha,2\beta)-\p_yF(0,2\beta) \end{array}\right).
\end{equation}
\end{Lem}
\begin{proof}
The symmetry of initial data \eqref{eq:intial data 4v} and the symmetry of \eqref{eq:CGLODE} imply that we can take the ansatz that the solution $\va$ satisfies \eqref{eq:symmetry of a 4v}. Substituting \eqref{eq:symmetry of a 4v} into \eqref{eq:CGLODE} and noting \eqref{eq:derivative of F on middle},  \eqref{eq:derivative of W}, we have
\begin{align}
  \left(\begin{array}{cc}1&-\lambda\\\lambda&1\end{array}\right)\dot{\va}_1&=2(\nabla F(2\alpha,2\beta)-\nabla F(2\alpha,0)-\nabla F(0,2\beta))\nn&=2\left(\begin{array}{c}\p_xF(2\alpha,2\beta)-\p_xF(2\alpha,0)\\\p_yF(2\alpha,2\beta)-\p_yF(0,2\beta)\end{array}\right),
\end{align}
which together with \eqref{eq:symmetry of a 4v} immediately implies \eqref{eq:sol of alpha beta 4v}.
\end{proof}
\begin{figure}[htp!]
\begin{center}
\begin{tabular}{cc}
\includegraphics[height=6cm]{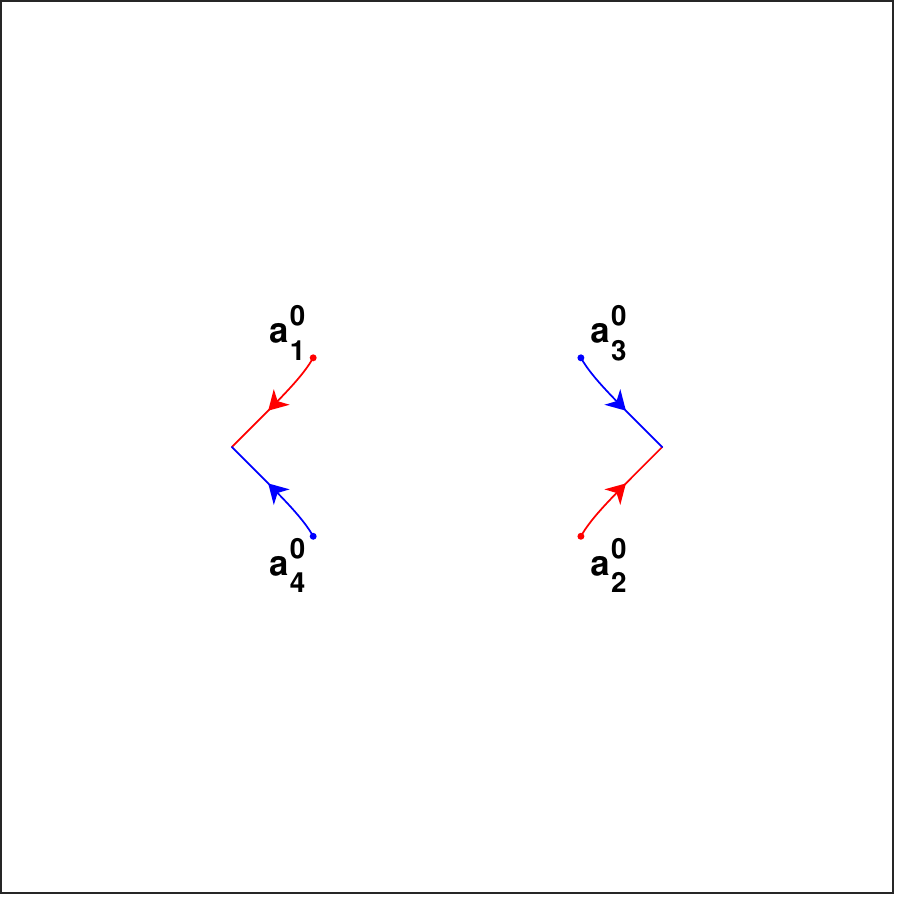}
&\includegraphics[width=6cm]{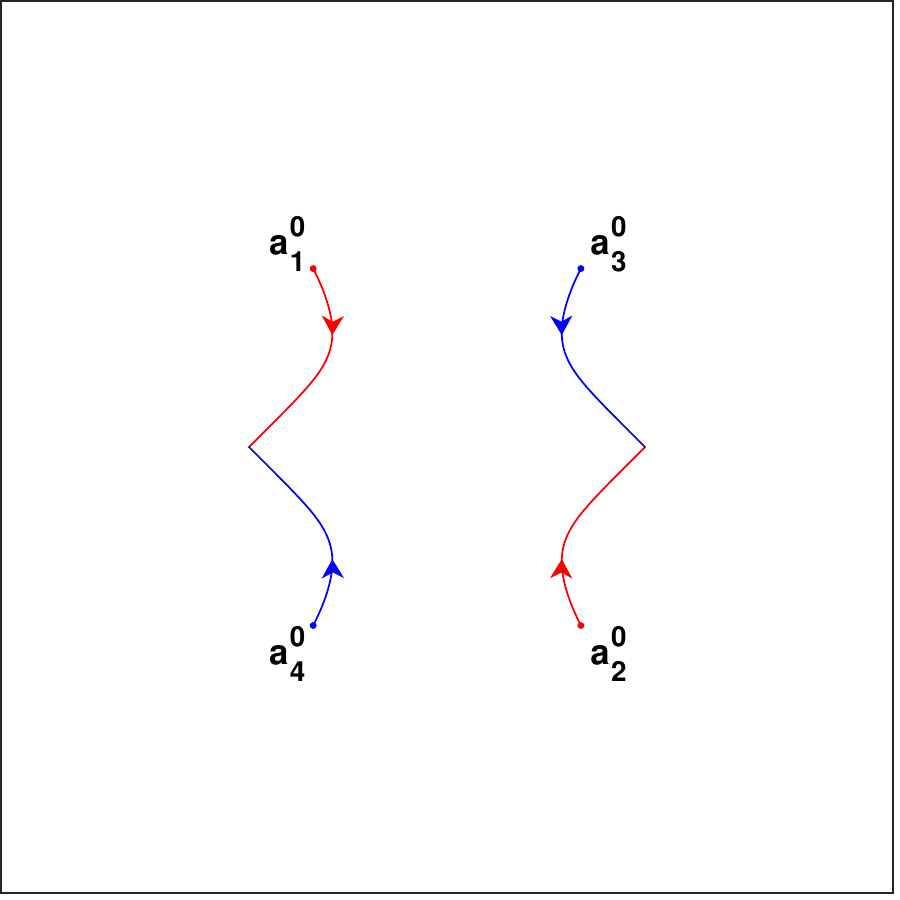}

\\
$(\alpha^0,\beta^0)=(-0.15,0.1)$&$(\alpha^0,\beta^0)=(-0.15,0.2)$\\
\includegraphics[height=6cm]{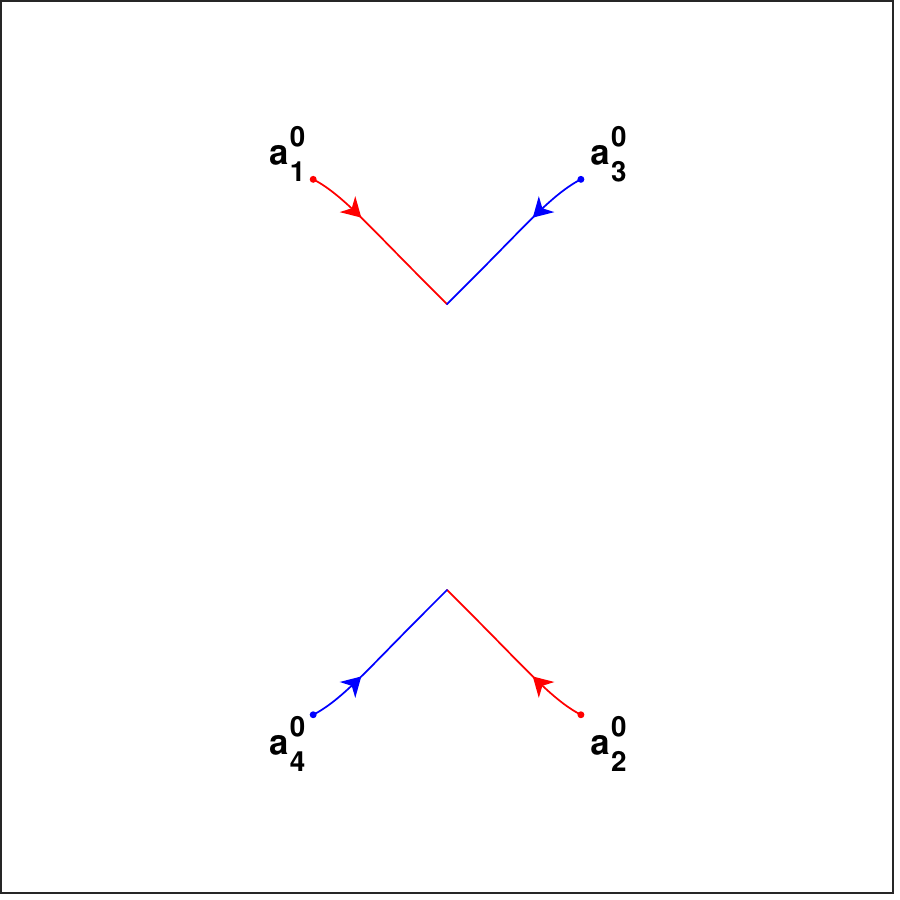}&
\includegraphics[height=6cm]{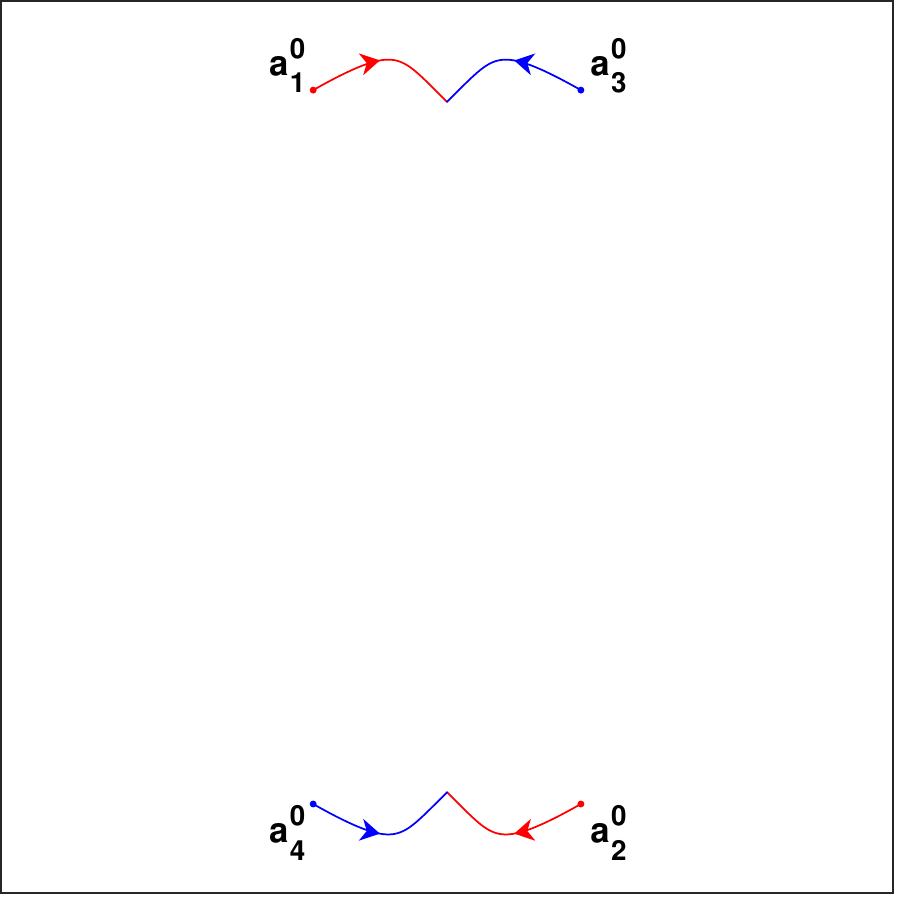}\\
$(\alpha^0,\beta^0)=(-0.15,0.3)$&$(\alpha^0,\beta^0)=(-0.15,0.4)$\\
\includegraphics[height=6cm]{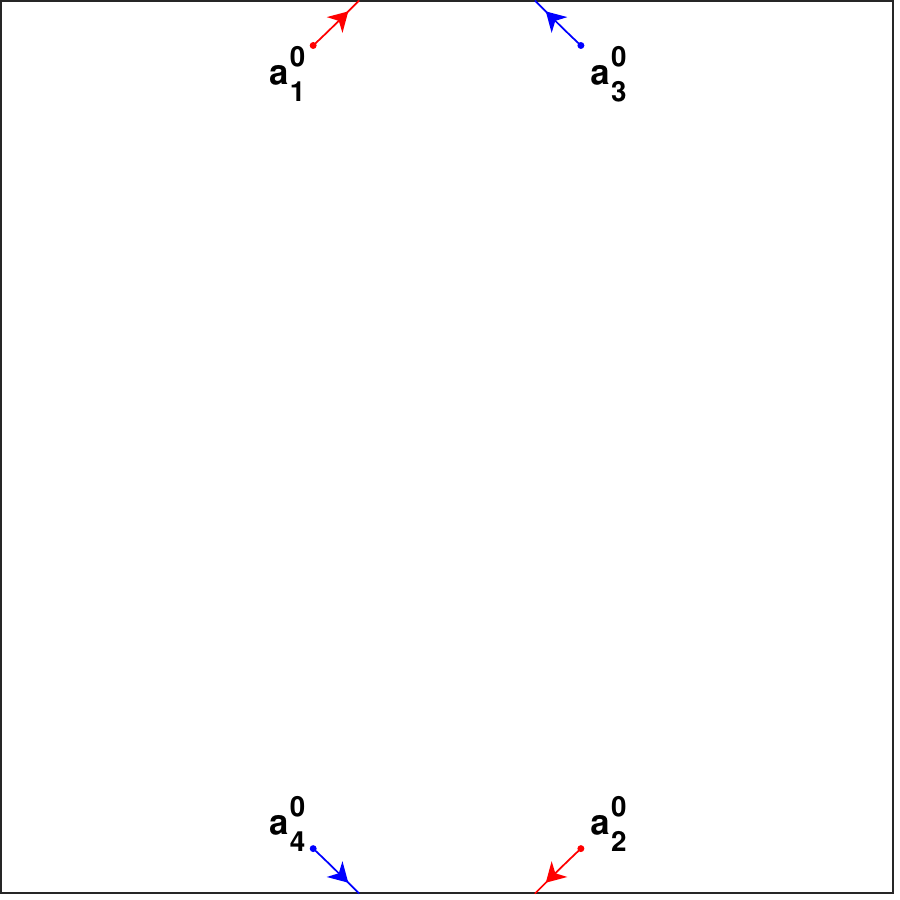}
&\includegraphics[height=6cm]{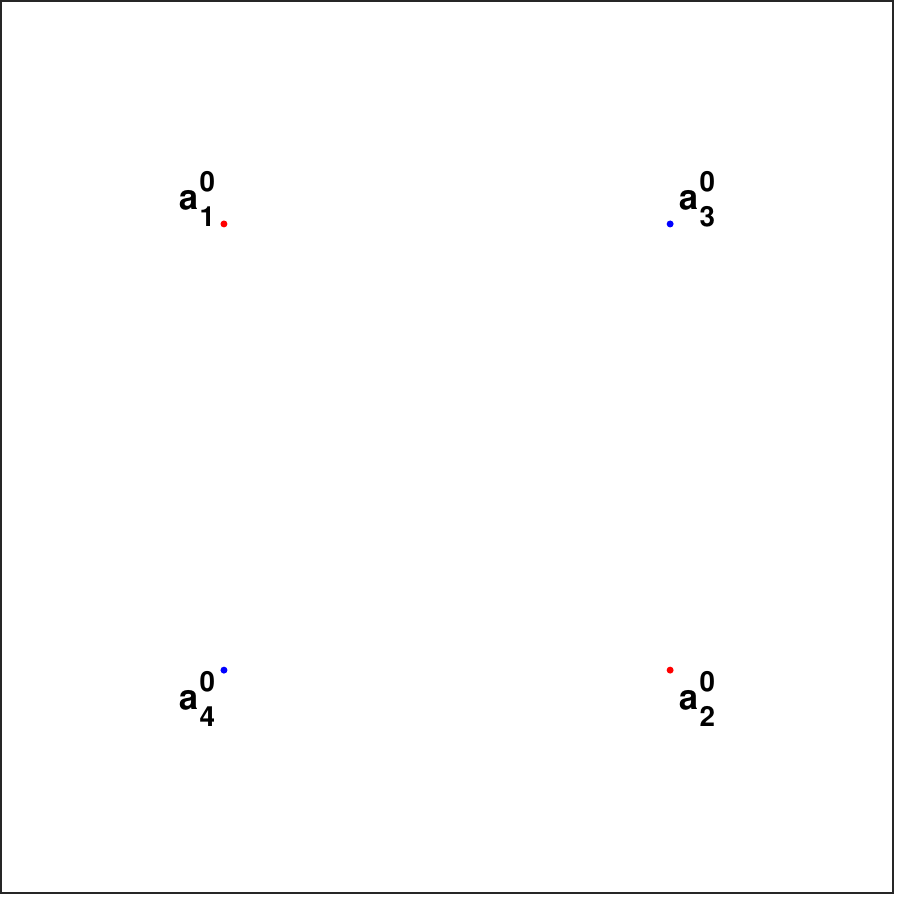}
\\
$(\alpha^0,\beta^0)=(-0.15,0.45)$&$(\alpha^0,\beta^0)=(-0.25,0.25)$
\end{tabular}
\end{center}
\caption{Some trajectories of solution of \eqref{eq:CGLODE} with $\lambda=1$ and initial data satisfying \eqref{eq:intial data 4v}. 
}\label{fig:path 4 v}
\end{figure}

Figure \ref{fig:path 4 v} shows some numerical results for the solution of \eqref{eq:CGLODE} with $\lambda=1$ and initial data satisfying \eqref{eq:symmetry of a}. The first five pictures in Figure \ref{fig:path 4 v} give five typical trajectories for fixed $\alpha^0$: (i) when $\beta^0$ is small, $\va_1$ will collide with $\va_4$ while $\va_2$ will collide with $\va_3$; (ii) as $\beta^0$ increases, $\va_1$ will get closer to $\va_3$, and if $\beta^0$ is large enough, $\va_1$ will collide with $\va_3$ while $\va_2$ will collide with $\va_4$; (iii) when $\beta^0$ is close to $0.5$ enough, $\va_1$ will collide with $\va_4$ again but in a direction different with the case (i). The last picture in Figure \ref{fig:path 4 v} gives an equilibrium state, which is caused by \eqref{eq:derivative of F on edge} and \eqref{eq:derivative of F on middle}.

\section{Conclusion}
The reduced dynamical laws for quantized vortex dynamics of the complex Ginzburg-Landau equation on torus were established when the core size of vortex $\varepsilon\to 0$. The motion of vortices is governed by a mixed flow of gradient flow and Hamiltonian flow which are both driven by a renormalized energy on torus. Finally, a first integral of the reduced dynamical laws was presented and some analytical solutions with several initial setups with symmetry were obtained.

\section*{CRediT authorship contribution statement}
\textbf{Yongxing Zhu:} Conceptualization, Methodology, Writing - original draft, Writing - review \& editing.

\section*{Data availability}
No data was used for the research described in the article.

\section*{Acknowledgments}
This work was partially supported by the China Scholarship Council (Y. Zhu). Part of the work was done when the author was visiting National University of Singapore during 2021-2023 and the Institute for Mathematical Science in 2023. The author would like to express his sincere gratitude to Prof. Huaiyu Jian in Tsinghua University and Prof. Weizhu Bao in National University of Singapore for their guidance and encouragement.

\section*{Declaration of competing interest}
The authors declare that they have no known competing financial interests or personal relationships that could have appeared to influence the work reported in this paper.

\end{document}